\documentclass{article} 
\usepackage{amssymb,latexsym,pstricks,pst-node,pst-coil}

\def\neweq{\setcounter{equation}{0}} 
\newtheorem{theorem}{Theorem}[section] 
\newtheorem{proposition}[theorem]{Proposition} 
\newtheorem{corollary}[theorem]{Corollary} 
\newtheorem{definition}[theorem]{Definition} 
 
\newtheorem{remark}[theorem]{Remark} 
\newtheorem{lemma}[theorem]{Lemma} 
 
\newtheorem{example}[theorem]{Example}

\def\Z{\mathbb{Z}}
\def\RR{\mathbb{R}}
\def\C{\mathcal{C}}
\def\L{\mathcal{L}}
\def\Id{\mathrm{Id}}
\def\R{\mathcal{R}}

\def\a{\mathbf{a}}

\def\addr{\footnotesize \it}
\def\b{\mathbf{b}}

\def\l{\ell}
\def\pprec{\cdot\kern-.2in\prec}
\def\rmleft{\mathrm{l}}
\def\rmright{\mathrm{r}}
\def\sct{\scriptstyle}
\def\wnot{w_\mathrm{o}}
\def\anot{{\mathbf{a}_\mathrm{o}}}
\def\({\left(}
\def\){\right)}

\def\proof{\smallskip\noindent {\it Proof --- \ }} 
\def\proofof#1{\smallskip\noindent {\it Proof of #1 --- \ }} 
\def\endproof{\hfill$\square$\medskip}

\title{Littlewood-Richardson Coefficients via Yang-Baxter Equation}

\author{%
\begin{tabular}{ccc}
Oleg Gleizer & \qquad & Alexander Postnikov  \\[.05in]
{\addr Northeastern University, Boston, MA 02115} & &
{\addr University of California, Berkeley, CA 94720} \\[.2in]
\end{tabular}}

\date{{\normalsize August~30, 1999\footnote{%
The latest version of this paper can be found on the web at
\texttt{http://www.math.berkeley.edu/\~{}apost/}}}}
 
\begin{document} 
\maketitle

\begin{abstract} 
The purpose of this paper is to present an interpretation for the decomposition
of the tensor product of two or more irreducible representations of~$GL(N)$ in
terms of a system of quantum particles.  Our approach is based on a certain
scattering matrix that satisfies a Yang-Baxter type equation.  The
corresponding piecewise-linear transformations of parameters 
give a solution to the tetrahedron equation. 
These transformation maps are naturally related to the dual canonical bases for
modules over the quantum enveloping algebra~$U_q(sl_n)$.  A byproduct of our
construction is an explicit description for the cone of Kashiwara's parametrizations of
dual canonical bases. This solves a problem posed by Berenstein and
Zelevinsky.  We present a graphical interpretation of the scattering matrices
in terms of web functions, which are related to honeycombs of Knutson and Tao.
\end{abstract}

\section{Introduction}
\label{sec:intro}
\neweq 

The aim of this paper is to further investigate the Grothendieck ring~$K_N$ of
polynomial  representations of the general linear group~$GL(N)$.  
Let~$V_\lambda$ be the irreducible representation of $GL(N)$
with highest weight~$\lambda$.  The
structure constants~$c_{\lambda\mu}^\nu$ of the Grothendieck ring in the basis
of irreducible representations are given by
$$
V_\lambda\otimes V_\mu= \sum_{\nu} c_{\lambda\mu}^\nu V_\nu\,.
$$
We also mention several alternative interpretations for the
numbers~$c_{\lambda\mu}^\nu$.  These numbers are:

\begin{itemize}
\item the structure constants of the ring of symmetric polynomials in the basis of
Schur polynomials;
\item the coefficients of the decomposition into irreducibles of representations
of symmetric groups induced from parabolic subgroups;
\item the structure constants of the cohomology ring of a Grassmannian in the
basis of Schubert classes.
\end{itemize}

The celebrated Littlewood-Richardson rule (see, e.g.,~\cite{macdonald}) is an
explicit combinatorial description of the coefficients~$c_{\lambda\mu}^\nu$.
Several variations of this rule are known, including Zelevinsky pictures and
Berenstein-Zelevinsky triangles~\cite{BZ1}.

\medskip

In this paper we present a new interpretation of the Grothendieck ring~$K_N$
and the Littlewood-Richardson coefficients~$c_{\lambda\mu}^\nu$.
Our construction is based on the {scattering matrix} $R(c)$ that acts
in the tensor square of the linear space~$E$ with the basis $e_0, e_1, e_2, \dots$ by
$$
R(c)\,:\, e_x\otimes e_y \longmapsto 
\left\{
\begin{array}{cl}
e_{y+c}\otimes e_{x-c} & \textrm{if } c \geq x-y, \\[.1in]
0                      & \textrm{otherwise.}
\end{array}
\right.
$$
(We assume that $e_x=0$ whenever $x<0$.)  
We denote by $R_{ij}(c)$ the operator acting on $E^{\otimes m}$ as $R(c)$ 
on the $i$-th and $j$-th copy of~$E$ and as an identity elsewhere.  
The tensor product of two irreducible
representations $V_\lambda$ and $V_\mu$ can be written as a certain combination
of the operators $R_{ij}(c)$.

Using the operators $R_{ij}(c)$ we define a new bilinear operation ``$*$'' 
on the tensor
algebra~$T(E)$ that will correspond to the operation of 
tensor product of representations of~$GL(N)$.  
It is straightforward that a Pieri type formula holds for the
$*$-product of any basis element in~$T(E)$ with~$e_k$.
The proof of the statement that tensor product $V_\lambda\otimes V_\mu$ is 
given by the $*$-product
easily follows from this fact and the fact that ``$*$'' is an associative 
operation.

The associativity of the $*$-multiplication
is obtained
from the following Yang-Baxter type relation for the scattering matrices.  
The operators $R_{12}(c_{12})$, $R_{13}(c_{13})$, and $R_{23}(c_{23})$ 
acting on~$E^{\otimes 3}$ satisfy the relation 
$$
R_{23}(c_{23})R_{13}(c_{13})R_{12}(c_{12}) =
R_{12}(c_{12}')R_{13}(c_{13}')R_{23}(c_{23}'),
$$
where~$c_{12}, c_{13}, c_{23}$ are arbitrary parameters and
$c_{12}', c_{13}', c_{23}'$ are given by the following piecewise-linear
formulas
\begin{equation}
\begin{array}{l}
c_{12}'=\min(c_{12},c_{13}-c_{23}),\\[.05in]
c_{13}'=c_{12}+c_{23},\\[.05in]
c_{23}'=\max(c_{23},c_{13}-c_{12}).
\end{array}
\label{eq:intro-transform}
\end{equation}

\medskip

Surprisingly, the same piecewise-linear transformations arise
in the theory of dual canonical bases for the modules over the 
quantum enveloping algebra~$U_q(sl_n)$ (see~\cite{BZ0, BZ2}).
For a fixed reduced decomposition of the longest element~$\wnot$
in the symmetric group~$S_n$, elements of the dual canonical 
basis (also known as the string basis) are parameterized by ${n\choose 2}$-tuples of
integers (strings) that belong to a certain string cone
(Kashiwara's parametrization).  Two
parametrizations that correspond to reduced decompositions related by a
Coxeter move 
are obtained from each other by the formulas~(\ref{eq:intro-transform}).

The string cone was described in~\cite{BZ0} for a certain reduced 
decomposition of~$\wnot$.  The core of our construction lies
in an explicit description of the string cone for any reduced 
decomposition.  Thus we solve a rather nontrivial problem 
posed in~\cite{BZ0}.

We also present a graphical (or ``pseudo-physical'') interpretation of the
scattering matrices and their compositions in the language of web functions
and ``systems of quantum particles.''
Web functions are closely related to honeycombs of
Knutson and Tao~\cite{Knutson} and
Berenstein-Zelevinsky triangles~\cite{BZ1}.
It is shown in~\cite{Knutson} that integral honeycombs are in one-to-one correspondence 
with Berenstein-Zelevinsky 
patterns.   We establish a simple ``dual'' correspondence between integral web functions
and Berenstein-Zelevinsky patterns.  This reveals the ``hidden duality'' of 
the Littlewood-Richardson coefficients under the conjugation of partitions.

\medskip

We briefly outline the structure of the paper.  In
Section~\ref{sec:preliminaries} we give some background on the representation
theory of general linear groups, the Littlewood-Richardson coefficients, and
the combinatorics of symmetric groups and reduced decompositions.  In
Section~\ref{sec:scattering} we define the scattering matrices~$R_{ij}(c)$ and
formulate our rule for the Littlewood-Richardson coefficients.
Section~\ref{sec:yb} is devoted to the Yang-Baxter type relation for
the scattering matrices.
In Section~\ref{sec:cone} we define and study principal cones for
any reduced decomposition of a permutation. 
In the case of the longest permutation these cones are exactly the 
string cones of parametrizations of dual canonical bases.
The associativity of the $*$-product is deduced in
Section~\ref{sec:associativity}.  In Section~\ref{sec:web}
we define web functions and establish their relationship with
the scattering matrices and Berenstein-Zelevinsky patterns.

\section{Preliminaries}                
\label{sec:preliminaries}
\neweq

In this section we remind the reader the basic notions and notation 
related to 
symmetric groups and representations of general linear groups.

\subsection{Representations of general linear groups}

Let us recall the basics of the representation theory of the
general linear group~$GL(N)$.

The general linear group~$GL(N)$ is the automorphism group of the
$N$-dimensional complex linear space~$\mathbb{C}^N$.  A complex
finite-dimensional linear representation~$V$ of~$GL(N)$ is called {\it
polynomial} if the corresponding mapping $GL(N)\to\mathrm{Aut}(V)$ is given by
polynomial functions.  An arbitrary holomorphic finite-dimensional 
representation is obtained by
tensoring a polynomial representation with a determinant representation
$g\mapsto\det^k(g)$ for suitable negative $k$.

An irreducible polynomial representation of~$GL(N)$ is uniquely determined
by its highest weight~$\lambda=(\lambda_1,\dots,\lambda_N)$,
which can be any integer element of the dominant 
chamber given by $\lambda_1\geq \lambda_2 \geq \cdots \geq\lambda_N\geq 0$.  
We denote by~$V_\lambda$ the irreducible representation with highest 
weight~$\lambda$.  Its degree is~$|\lambda|=\lambda_1+\cdots+\lambda_N$.

The collection of polynomial representations of~$GL(N)$ equipped with the
operations of direct sum and tensor product has the structure of an abelian
category.  Let~$K_N=K(GL(N))$ be the Grothendieck ring of this category.
Degree of representations provides a natural grading on the ring~$K_N$.
Slightly abusing notations, we will identify a representation with its
image in the Grothendieck ring~$K_N$.

The irreducible representations~$V_\lambda$ form a $\Z$-basis of~$K_N$.  Our
primary interest is in the structure constants of~$K_N$.  In other words, we
would like to investigate the coefficients~$c_{\lambda\mu}^\nu$ of the
expansion of the tensor product of two irreducible representations into a 
direct sum of irreducibles:
$$
V_\lambda\otimes V_\mu= \sum_{\nu} c_{\lambda\mu}^\nu V_\nu\,.
$$

The weights $\omega_k=(1,\dots,1,0,\dots,0)$ (with $k$ ones) are called the
{\it fundamental weights.}  By convention $\omega_0=(0,\dots,0)$.  Every
dominant weight~$\lambda$ can be written uniquely as a sum of fundamental
weights $\lambda=\omega_{x_1}+\cdots+\omega_{x_m}$, $1\leq x_1\leq \dots\leq
x_m\leq N$.  
Actually, the numbers $x_i$ are just parts of the partition $\lambda'$ 
{\it conjugate} to~$\lambda$, i.e., $\lambda'=(x_m,x_{m-1},\dots,x_1)$.

The fundamental representation~$V_{\omega_k}$ is the $k$-th exterior power of
the tautological representation of $GL(N)$.  
Pieri's formula gives an explicit rule for the tensor product
of~$V_{\omega_k}$ with an irreducible representation~$V_\lambda$.  

\begin{proposition} {\rm (Pieri's formula) } \ 
For $\lambda=\omega_{x_1}+\cdots+\omega_{x_m}$, $1\leq x_1\leq x_2\leq \cdots
\leq x_m\leq N$,
we have
\begin{equation} 
V_{\omega_k}\otimes V_\lambda = \sum V_\mu\,,
\label{eq:pieri}
\end{equation} 
where the sum is over all $\mu=\omega_{y_1}+\cdots +\omega_{y_{m+1}}$ satisfying
the following interlacing conditions:
$$
\begin{array}{l}
0\leq y_1\leq x_1\leq y_2\leq x_2\cdots\leq y_m\leq x_m\leq y_{m+1}\leq N,\\[.1in]
y_1-x_1+y_2-x_2+\cdots+y_m-x_m+y_{m+1}=k.
\end{array}
$$
\end{proposition}

The Grothendieck ring~$K_N$ is generated by the fundamental representations 
$V_{\omega_k}$.  This implies the following statement that will be handy afterward.

\begin{lemma}
\label{le:odot}
Suppose that $\odot$ is a bilinear associative multiplication operation
on the linear space $K_N$ such that for any fundamental weight $\omega_k$ 
and any dominant weight $\lambda$ the product $V_{\omega_k}\odot V_\lambda$
is given by  Pieri's formula~{\rm (\ref{eq:pieri})} and $V_{(0,\dots,0)}$
is the identity element.  Then $\odot$ is the usual
multiplication in~$K_N$---tensor product of representations.
\end{lemma}

\proof  
We will show that $V_\lambda\odot V_\mu=V_\lambda\otimes V_\mu$
by induction on the degree $|\lambda|$ of~$V_\lambda$.
First, $V_{(0,\dots,0)}\odot V_\mu = V_\mu$ by the condition of lemma.  
Suppose that the statement is true for any~$V_\lambda$ with $|\lambda|<d$.
For $|\lambda|=d$, we can express $V_\lambda$ via the generators~$V_{\omega_k}$ 
as $\sum V_{\omega_k}\otimes W_k$ where the~$W_k$ are
degree~$d-1$ elements of~$K_N$. Then, by the inductive hypothesis, 
$$
V_\lambda\odot V_\mu= (\sum V_{\omega_k}\odot W_k)\odot V_\mu=
\sum V_{\omega_k}\odot (W_k\odot V_\mu) = \sum V_{\omega_k}\otimes
W_k\otimes V_\mu=V_\lambda\otimes V_\mu.
$$
\endproof

\subsection{Symmetric group}

Our constructions rely strongly on the combinatorics of reduced decompositions
in the symmetric group~$S_n$.  This section is devoted to a brief account
of this theory.

Let $s_a\in S_n$ be the adjacent transposition that interchanges $a$ and $a+1$.
Then~$s_1,\dots,s_{n-1}$ generate the symmetric group~$S_n$.   The
generators~$s_a$ satisfy the following Coxeter relations: 
\begin{equation} 
\begin{array}{l}
s_a^2=1,\\[.1in]
s_a s_b = s_b s_a,\quad\textrm{for }|a-b|\geq 2,\\[.1in]
s_a s_{a+1} s_a=s_{a+1} s_a s_{a+1}.
\end{array}
\label{eq:coxeter}
\end{equation} 

For a permutation $w\in S_n$, an expression $w=s_{a_1}s_{a_2}\cdots s_{a_l}$
of minimal possible length~$l$ is called a {\it reduced decomposition};
and $l=\l(w)$ is the {\it length} of~$w$.  The corresponding sequence
$\a=(a_1,a_2,\dots,a_l)$ is called a {\it reduced word} for~$w$.  Let $\R(w)$
denote the set of all reduced words for~$w$.  A pair $(i,j)$, $1\leq i<j\leq
m$, is called an {\it inversion} in~$w$ if $w(i)>w(j)$. By $I(w)$ we denote
the set of all inversions of~$w$.  The number~$|I(w)|$ of inversions in~$w$ is
equal to its length~$\l(w)$.  

Let $\wnot$ be the longest permutation in~$S_n$ given by $\wnot(i)=n+1-i$.
Then $I(\wnot)$ is the set of all pairs $1\leq i<j\leq n$.
A total ordering ``$<$'' of inversions $(i,j)$  in $I(\wnot)$ is said to be a 
{\it reflection ordering} if for any triple $i<j<k$ we have 
$$
\textrm{either}\quad 
(i,j)< (i,k)< (j,k)
\quad\textrm{or}\quad
(j,k)< (i,k)< (i,j).
$$
Also, for any~$w\in S_n$, we say that a total ordering of inversions in $I(w)$
is a {\it reflection ordering} if it is a final interval of some reflection
ordering of $I(\wnot)$.

The set of all reflection orderings of~$I(w)$ is in one-to-one correspondence
with the set of reduced decompositions of~$w$, cf.~\cite[Proposition~2.13]{Dy}.
Explicitly, for a reduced decomposition $w=s_{a_1}s_{a_2}\cdots s_{a_l}$, the
sequence of pairs $(i_1,j_1)< \cdots < (i_l,j_l)$ such that $i_r=s_{a_l}
s_{a_{l-1}}\cdots s_{a_{r+1}}(a_r)$ and $j_r=s_{a_l} s_{a_{l-1}}\cdots
s_{a_{r+1}}(a_r+1)$, $r=1,\dots,l$, is a reflection ordering of~$I(w)$.
Moreover, every reflection ordering of~$I(w)$ arises in this fashion.

Graphically, we represent a reduced decomposition by its {\it wiring diagram,}
which is also called a {\it pseudo-line arrangement.}
For instance, the reduced decomposition $s_3\, s_2\, s_1\, s_2$ 
of an element in~$S_4$ is depicted by the diagram
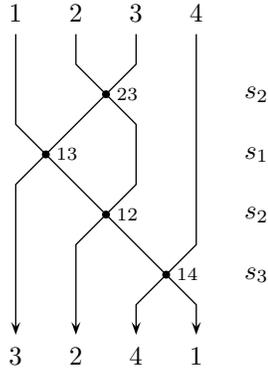
\psset{curvature=0 0 0}
\psset{xunit=.4cm,yunit=.4cm}
\begin{center}
\pspicture[.1](0,0)(10,12)

  \pscurve[linewidth=.5pt, showpoints=false]{<-}%
          (2,1)(2,2)(2,3)(2,4)(2,5)(2,6)(3,7)(4,8)(5,9)(6,10)(6,11)
  \pscurve[linewidth=.5pt, showpoints=false]{<-}%
          (4,1)(4,2)(4,3)(4,4)(5,5)(6,6)(6,7)(6,8)(5,9)(4,10)(4,11)
  \pscurve[linewidth=.5pt, showpoints=false]{<-}%
          (6,1)(6,2)(7,3)(8,4)(8,5)(8,6)(8,7)(8,8)(8,9)(8,10)(8,11)
  \pscurve[linewidth=.5pt, showpoints=false]{<-}%
          (8,1)(8,2)(7,3)(6,4)(5,5)(4,6)(3,7)(2,8)(2,9)(2,10)(2,11)
  \qdisk(7,3){1.5pt}
  \qdisk(5,5){1.5pt}
  \qdisk(3,7){1.5pt}
  \qdisk(5,9){1.5pt}
  \rput[u](2,11.7){$1$}
  \rput[u](4,11.7){$2$}
  \rput[u](6,11.7){$3$}
  \rput[u](8,11.7){$4$}
  \rput[u](2,.3){$3$}
  \rput[u](4,.3){$2$}
  \rput[u](6,.3){$4$}
  \rput[u](8,.3){$1$}
  \rput[u](5.7,9){$\sct 23$}
  \rput[u](3.7,7){$\sct 13$}
  \rput[u](5.7,5){$\sct 12$}
  \rput[u](7.7,3){$\sct 14$}
  \rput[u](10,9){$s_2$}
  \rput[u](10,7){$s_1$}
  \rput[u](10,5){$s_2$}
  \rput[u](10,3){$s_3$}
\endpspicture
\end{center}
The nodes of this diagram correspond to the adjacent transpositions.  On the
other hand, each node is a crossing of $i$-th and $j$-th pseudo-lines, where
$(i,j)$ forms an inversion.   Reading these pairs in the wiring diagram
from bottom to top gives the corresponding reflection ordering of the inversions.
In the above example, the associated reflection ordering is 
$(1,4)< (1,2)< (1,3)< (2,3)$.

Applying the Coxeter relations to reduced decompositions results in the local
transformations that are called $2$-moves and $3$-moves.  Namely, $2$-moves
correspond to the second equation in~(\ref{eq:coxeter}) and $3$-moves to the
third equation in~(\ref{eq:coxeter}).  Two reduced decompositions of the same
permutation are always connected by a sequence of $2$- and $3$-moves.
Graphically,  $2$- and $3$-moves can be represented by the following local
transformations of wiring diagrams, where $i<j<k<l$.

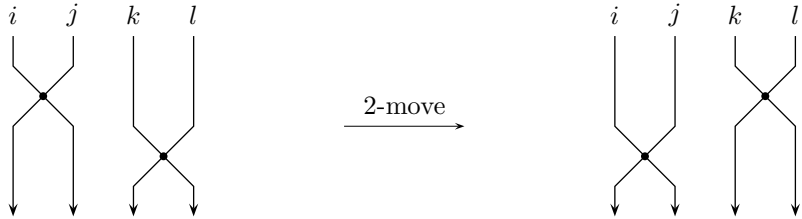
\psset{curvature=0 0 0}
\psset{xunit=.4cm,yunit=.4cm}
\begin{center}
\pspicture[.1](0,0)(29,8)

  \pscurve[linewidth=.5pt, showpoints=false]{<-}%
          (2,1)(2,2)(2,3)(2,4)(3,5)(4,6)(4,7)
  \pscurve[linewidth=.5pt, showpoints=false]{<-}%
          (4,1)(4,2)(4,3)(4,4)(3,5)(2,6)(2,7)
  \pscurve[linewidth=.5pt, showpoints=false]{<-}%
          (6,1)(6,2)(7,3)(8,4)(8,5)(8,6)(8,7)
  \pscurve[linewidth=.5pt, showpoints=false]{<-}%
          (8,1)(8,2)(7,3)(6,4)(6,5)(6,6)(6,7)
  \qdisk(3,5){1.5pt}
  \qdisk(7,3){1.5pt}
  \rput[u](2,7.7){$i$}
  \rput[u](4,7.7){$j$}
  \rput[u](6,7.7){$k$}
  \rput[u](8,7.7){$l$}

  \psline[linewidth=.2pt]{->}(13,4)(17,4)
  \rput[u](15,4.7){$2$-move}

  \pscurve[linewidth=.5pt, showpoints=false]{<-}%
          (22,1)(22,2)(23,3)(24,4)(24,5)(24,6)(24,7)
  \pscurve[linewidth=.5pt, showpoints=false]{<-}%
          (24,1)(24,2)(23,3)(22,4)(22,5)(22,6)(22,7)
  \pscurve[linewidth=.5pt, showpoints=false]{<-}%
          (26,1)(26,2)(26,3)(26,4)(27,5)(28,6)(28,7)
  \pscurve[linewidth=.5pt, showpoints=false]{<-}%
          (28,1)(28,2)(28,3)(28,4)(27,5)(26,6)(26,7)
  \qdisk(23,3){1.5pt}
  \qdisk(27,5){1.5pt}
  \rput[u](22,7.7){$i$}
  \rput[u](24,7.7){$j$}
  \rput[u](26,7.7){$k$}
  \rput[u](28,7.7){$l$}

\endpspicture
\end{center}
and

\psset{curvature=0 0 0}
\psset{xunit=.4cm,yunit=.4cm}
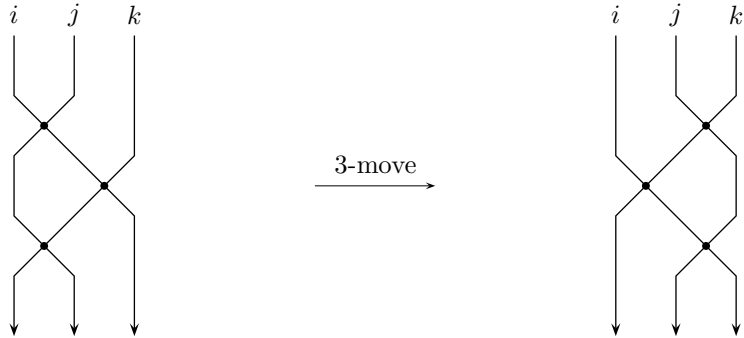
\begin{center}
\pspicture[.1](0,0)(27,12)

  \pscurve[linewidth=.5pt, showpoints=false]{<-}%
          (2,1)(2,2)(2,3)(3,4)(4,5)(5,6)(6,7)(6,8)(6,9)(6,10)(6,11)
  \pscurve[linewidth=.5pt, showpoints=false]{<-}%
          (4,1)(4,2)(4,3)(3,4)(2,5)(2,6)(2,7)(3,8)(4,9)(4,10)(4,11)
  \pscurve[linewidth=.5pt, showpoints=false]{<-}%
          (6,1)(6,2)(6,3)(6,4)(6,5)(5,6)(4,7)(3,8)(2,9)(2,10)(2,11)
  \qdisk(3,4){1.5pt}
  \qdisk(5,6){1.5pt}
  \qdisk(3,8){1.5pt}
  \rput[u](2,11.7){$i$}
  \rput[u](4,11.7){$j$}
  \rput[u](6,11.7){$k$}

  \psline[linewidth=.2pt]{->}(12,6)(16,6)
  \rput[u](14,6.7){$3$-move}

  \pscurve[linewidth=.5pt, showpoints=false]{<-}%
          (22,1)(22,2)(22,3)(22,4)(22,5)(23,6)(24,7)(25,8)(26,9)(26,10)(26,11)
  \pscurve[linewidth=.5pt, showpoints=false]{<-}%
          (24,1)(24,2)(24,3)(25,4)(26,5)(26,6)(26,7)(25,8)(24,9)(24,10)(24,11)
  \pscurve[linewidth=.5pt, showpoints=false]{<-}%
          (26,1)(26,2)(26,3)(25,4)(24,5)(23,6)(22,7)(22,8)(22,9)(22,10)(22,11)
  \qdisk(25,4){1.5pt}
  \qdisk(23,6){1.5pt}
  \qdisk(25,8){1.5pt}
  \rput[u](22,11.7){$i$}
  \rput[u](24,11.7){$j$}
  \rput[u](26,11.7){$k$}
\endpspicture
\end{center}

\section{Scattering Matrix}             
\label{sec:scattering}
\neweq

Let~$E$ be the linear space with a basis $e_x$, $x\in\Z_+$.  We will
always assume that $e_x=0$ for $x<0$.

\begin{definition} {\rm
For $c\in\Z$, the {\it scattering matrix} $R(c)$
is the linear operator which acts on the space $E\otimes E$ by
\begin{equation}
R(c)\,:\, e_x\otimes e_y \longmapsto 
\left\{
\begin{array}{cl}
e_{y+c}\otimes e_{x-c} & \textrm{if } c \geq x-y, \\[.1in]
0                      & \textrm{otherwise.}
\end{array}
\right.
\label{eq:scattering}
\end{equation}
}
\label{def:scattering}
\end{definition}

The space~$E$ can be viewed as the space of states of a certain quantum
particle.  The basis vector~$e_x$ corresponds to a particle with energy
level~$x$.  
We will think of the scattering matrix~$R(c)$  as the result of the interaction 
of two particles with energy levels~$x$ and~$y$.  
Pictorially, we can represent it by the following ``Feinman diagram''

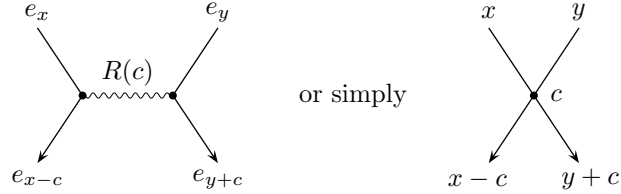
\psset{xunit=.3cm,yunit=.3cm}
\begin{center}
\pspicture[.1](0,0)(28,8)
  \psline[linewidth=.5pt]{<-}(2,1)(4,4)(2,7)
  \psline[linewidth=.5pt]{<-}(10,1)(8,4)(10,7)
  \pscoil[coilaspect=0,coilarm=0pt,coilheight=2, 
          coilwidth=2pt,linewidth=.3pt]{-}(4,4)(8,4)
  \qdisk(4,4){1.5pt}
  \qdisk(8,4){1.5pt}
  \rput[u](2,.3){$e_{x-c}$}
  \rput[u](10,.3){$e_{y+c}$}
  \rput[u](2,7.7){$e_x$}
  \rput[u](10,7.7){$e_y$}
  \rput[u](6,5){$R(c)$}
   
  \rput[u](16,4){or simply}

  \psline[linewidth=.5pt]{<-}(22,1)(26,7)
  \psline[linewidth=.5pt]{<-}(26,1)(22,7)
  \qdisk(24,4){1.5pt}
  \rput[u](21.5,.3){$x-c$}
  \rput[u](26.5,.3){$y+c$}
  \rput[u](22,7.7){$x$}
  \rput[u](26,7.7){$y$}
  \rput[u](25,4){$c$}
\endpspicture
\end{center}

Notice that the energy conservation law holds in our model, since the sum of
energies of particles after the interaction $(y+c)+(x-c)$ is the same as before
the interaction.

By $R_{ij}(c)$ we denote the linear endomorphism of $E^{\otimes m} =
E\otimes\cdots\otimes E$ which acts
as the transformation $R(c)$ on the $i$-th and the $j$-th copies of $E$ and 
as an identity operator on other copies.  
Let $\a=(a_1,\dots,a_l)$ be a reduced word for~$w\in S_n$, which
is associated with the reduced decomposition $w=s_{a_1}s_{a_2}\cdots s_{a_l}$,
and let $(i_1,j_1)< \cdots < (i_l,j_l)$ be the corresponding reflection 
ordering of the inversion set~$I(w)$.
For a collection~$C=(c_{ij})$, $(i,j)\in I(w)$,
of integer parameters, we define an endomorphism~$R_\a(C)$ of~$E^{\otimes m}$ 
as the composition of scattering matrices
\begin{equation}
\label{eq:rac}
R_\a(C) = R_{i_1j_1}(c_{i_1j_1})
R_{i_2j_2}(c_{i_2j_2})\cdots R_{i_lj_l}(c_{i_lj_l}).
\end{equation}
It is clear that $R_{ij}(c_{ij})$ commutes with $R_{kl}(c_{kl})$ 
provided that all $i,j,k,l$ are distinct.
Thus the composition~$R_\a(C)$ stays invariant when we apply
a $2$-move to the reduced word~$\a$.

For positive integers $m$ and~$n$, let~$w(m,n)$ be the permutation
from~$S_{m+n}$ given by
$$
\( \begin{array}{cccccccc}
1 & 2 & \cdots & m & m+1 & m+2 & \cdots & m+n \\
n+1 & n+2 & \cdots & n+m & 1 & 2 & \cdots & n 
\end{array}\).
$$
All reduced decompositions of the permutation~$w(m,n)$ are related by 
$2$-moves (cf.~the diagram below).  Thus the
map~$R_\a(C)$ does not depend upon any particular choice of a reduced
word~$\a$ for~$w(m,n)$.  We denote by $R_{(m,n)}(C)$ this endomorphism of
$E^{\otimes m} \otimes E^{\otimes n}$.  It depends upon the collection of
$mn$ parameters~$C=(c_{ij})$, $1\leq i\leq m<j\leq m+n$.

Let~$T(E)$ denote the tensor algebra of the linear space~$E$.   
We define a new bilinear operation $M:T(E)\otimes T(E)\to T(E)$ 
whose restriction $M_{m,n}:E^{\otimes m} \otimes E^{\otimes n} \to 
E^{\otimes (m+n)}$ 
is given by
\begin{equation}
M_{m,n} = \sum_C R_{(m,n)}(C),
\label{eq:mmn}
\end{equation} 
where the sum is over all collections $C$ of 
nonnegative integer parameters~$c_{ij}$, $1\leq i\leq m<j\leq m+n$, 
such that 
\begin{equation}
c_{ij}\geq c_{kl}\qquad\textrm{whenever}\quad k\leq i<j\leq l.
\label{eq:conditions}
\end{equation}
We will use the notation $A*B$ for $M(A,B)$, where $A,B\in T(E)$,
and occasionally call this multiplication operation {\it $*$-product.}  
Although the sum in~(\ref{eq:mmn}) involves infinitely many terms,
only a finite number of them are nonzero in the expansion for~$A*B$.

Let us remark that a collection $C$ of nonnegative integers
that satisfy~(\ref{eq:conditions}) is usually called  a {\it rectangular
shaped plane partition.}

The composition of scattering matrices $R_{(m,n)}$ can be represented by
the wiring diagram shown below (for $m=4$ and $n=3$).  The summation
in~(\ref{eq:mmn}) is over all collections of nonnegative integer 
parameters~$c_{ij}$ that weakly decrease downwards along the pseudo-lines 
of this diagram.

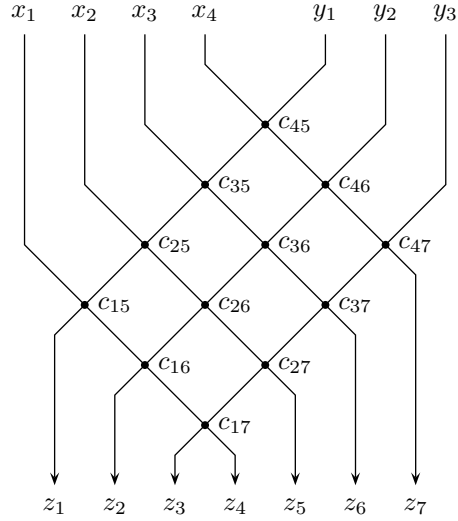
\psset{curvature=0 0 0}
\psset{xunit=.4cm,yunit=.4cm}
\begin{center}
\pspicture[.1](0,0)(16,17)

  \pscurve[linewidth=.5pt, showpoints=false]{<-}%
          (2,1)(2,6)(11,15)(11,16)
  \pscurve[linewidth=.5pt, showpoints=false]{<-}%
          (4,1)(4,4)(13,13)(13,16)
  \pscurve[linewidth=.5pt, showpoints=false]{<-}%
          (6,1)(6,2)(15,11)(15,16)

  \pscurve[linewidth=.5pt, showpoints=false]{<-}%
          (8,1)(8,2)(1,9)(1,16)
  \pscurve[linewidth=.5pt, showpoints=false]{<-}%
          (10,1)(10,4)(3,11)(3,16)
  \pscurve[linewidth=.5pt, showpoints=false]{<-}%
          (12,1)(12,6)(5,13)(5,16)
  \pscurve[linewidth=.5pt, showpoints=false]{<-}%
          (14,1)(14,8)(7,15)(7,16)

  \qdisk(7,3){1.5pt}
  \qdisk(9,5){1.5pt}
  \qdisk(11,7){1.5pt}
  \qdisk(13,9){1.5pt}
  \qdisk(5,5){1.5pt}
  \qdisk(7,7){1.5pt}
  \qdisk(9,9){1.5pt}
  \qdisk(11,11){1.5pt}
  \qdisk(3,7){1.5pt}
  \qdisk(5,9){1.5pt}
  \qdisk(7,11){1.5pt}
  \qdisk(9,13){1.5pt}

  \rput[u](8,3){$c_{17}$}
  \rput[u](6,5){$c_{16}$}
  \rput[u](4,7){$c_{15}$}
  \rput[u](10,5){$c_{27}$}
  \rput[u](8,7){$c_{26}$}
  \rput[u](6,9){$c_{25}$}
  \rput[u](12,7){$c_{37}$}
  \rput[u](10,9){$c_{36}$}
  \rput[u](8,11){$c_{35}$}
  \rput[u](14,9){$c_{47}$}
  \rput[u](12,11){$c_{46}$}
  \rput[u](10,13){$c_{45}$}

  \rput[u](1,16.7){$x_1$}
  \rput[u](3,16.7){$x_2$}
  \rput[u](5,16.7){$x_3$}
  \rput[u](7,16.7){$x_4$}

  \rput[u](11,16.7){$y_1$}
  \rput[u](13,16.7){$y_2$}
  \rput[u](15,16.7){$y_3$}

  \rput[u](2,.3){$z_1$}
  \rput[u](4,.3){$z_2$}
  \rput[u](6,.3){$z_3$}
  \rput[u](8,.3){$z_4$}
  \rput[u](10,.3){$z_5$}
  \rput[u](12,.3){$z_6$}
  \rput[u](14,.3){$z_7$}
\endpspicture
\end{center}

\begin{theorem}
\label{th:associative}
The space~$T(E)$ equipped with the multiplication 
operation $M$ is an associative ring.  
\end{theorem}

Recall that $\omega_1,\dots,
\omega_N$ are the fundamental weights of $GL(N)$.  By convention $\omega_0=0$.

\begin{theorem}  
\label{th:main}
The projection $p_N:T(E)\to K_N$ 
defined on the basis elements by
$$
p_N: e_{x_1}\otimes \cdots\otimes e_{x_m}\longmapsto
\left\{
\begin{array}{ll}
V_\lambda,\ \lambda=\omega_{x_1}+\cdots+\omega_{x_m}, &
            \textrm{provided }
            x_1\leq x_2 \leq \cdots \leq x_m\leq N,\\[.1in]
0         & \textrm{otherwise}
\end{array}
\right.
$$
is a homomorphism from the ring $(T(E), M)$ to the Grothendieck ring~$K_N$ of 
polynomial representations
of $GL(N)$.  In other words, if $p_N(A)=V_\lambda$ and $p_N(B)=V_\mu$ then
$p_N(A*B)=V_\lambda\otimes V_\mu$, the tensor product of representations.
\end{theorem}

Summarizing the above assertions and definitions, we can
formulate a rule for the Littlewood-Richardson coefficients.  
Let us denote by
$e_{x_1\dots x_m}$ the element $e_{x_1}\otimes \cdots \otimes e_{x_m} \in
T(E)$.  

\begin{corollary}
\label{cor:rule}
Let $\lambda = \omega_{x_1}+\cdots +\omega_{x_m}$,
$\mu = \omega_{y_1}+\cdots +\omega_{y_n}$, and
$\nu = \omega_{z_1}+\cdots +\omega_{z_{m+n}}$, where 
$x_1\leq \cdots \leq x_m$, 
$y_1\leq \cdots \leq y_n$, and
$z_1\leq \cdots \leq z_{m+n}$. 
The Littlewood-Richardson coefficient $c_{\lambda\mu}^\nu$ 
is equal to the number of collections $C$ of nonnegative integers
$c_{ij}$, $1\leq i\leq m$, $m+1\leq j\leq m+n$, such that 
$$
c_{ij}\geq c_{kl},\qquad \textrm{for } i\leq k <l \leq j;
$$
and 
$$
R_{(m,n)}(C) \cdot (e_{x_1\dots x_m}\otimes 
e_{y_1\dots y_n}) = e_{z_1\dots z_{m+n}}.
$$
\end{corollary}

\begin{proposition} 
\begin{enumerate}
\item
We have $e_{x_1\dots x_m}*e_{y_1\dots y_n}=0 $
unless $x_1\leq \cdots\leq x_m$ and $y_1\leq \cdots\leq y_n$.
\item
The product $e_{x_1\dots x_m}*e_{y_1\dots y_n}$ involves 
only terms $e_{z_1\dots z_{m+n}}$ with $z_1\leq \cdots\leq  z_{m+n}$.
\end{enumerate}
\end{proposition}

\proof
1. \ 
First, we show that applying $R_{13}(c_{13})R_{23}(c_{23})$ to 
$e_{x_1}\otimes e_{x_2}\otimes e_{y_1}$ always results in zero 
provided $c_{23}\geq c_{13}$ and $x_1>x_2$.  
Indeed, we have $R_{13}(c_{13}) R_{23}(c_{23})
\cdot (e_{x_1}\otimes e_{x_2}\otimes e_{y_1})=
R_{13}(c_{13})\cdot (e_{x_1}\otimes e_{y_1+c_{23}}\otimes e_{x_2-c_{23}})$ 
(or zero).  This expression is nonzero only if 
$c_{13}\geq x_1-(x_2-c_{23})$, i.e., $c_{13}-c_{23}\geq x_1-x_2$.  
Contradiction.

In general, suppose that, say, $x_i>x_{i+1}$.  The composition of operators
$R_{(m,n)}(C)$ with $C$ satisfying~(\ref{eq:conditions}) involves
the fragment $R_{i+1\,m+1}(c_{i+1\,m+1})R_{i\,m+1}(c_{i\,m+1})$,
where $c_{i+1\,m+1}\geq c_{i\,m+1}$.  By the above argument, applying 
these operators gives zero.

\medskip
\noindent
2. \ This statement follows by induction on~$m$ from
the next Proposition~\ref{pr:e-pieri}.
\endproof

Let us verify the statement of Theorem~\ref{th:main} for the $*$-product of
$e_x$ with an arbitrary $e_{x_1\dots x_m}$.  This product is given by the
following Pieri-type formula.

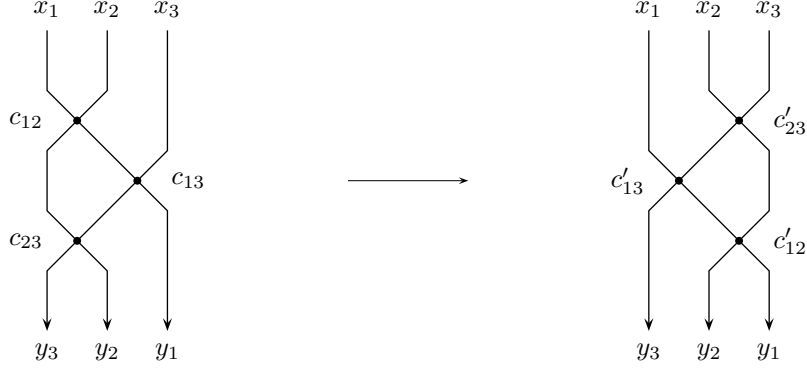
\begin{proposition}
\label{pr:e-pieri}
For $0\leq x_1\leq\cdots \leq x_m$, we have
$$
e_x*e_{x_1\dots x_m} = \sum e_{y_1\dots y_{m+1}},
$$
where the sum is over all $y_1,\dots, y_{m+1}$ satisfying
the following interlacing conditions:
\begin{equation}
\begin{array}{l}
0\leq y_1\leq x_1\leq y_2\leq x_2\cdots\leq y_m\leq x_m\leq y_{m+1},\\[.1in]
y_1-x_1+y_2-x_2+\cdots+y_m-x_m+y_{m+1}=x.
\end{array}
\label{eq:interlace}
\end{equation}
\end{proposition}

\proof
By definition, $e_x*e_{x_1\dots x_m} = 
\sum R_{1\,m+1}(c_m)R_{1\,m}(c_{m-1})\cdots R_{12}(c_1)\cdot (e_x\otimes
e_{x_1}\otimes\dots e_{x_m})$, where the sum is over
$c_1\geq c_2\geq\dots\geq c_m\geq 0$.  Each nonvanishing summand in the 
previous sum is equal to $e_{x-c_1}\otimes e_{x_1+c_1-c_2}\otimes
e_{x_2+c_2-c_3}\otimes\dots e_{x_{m-1}+c_{m-1}-c_m}\otimes e_{x_m+c_m}$
provided $c_1\geq x-x_1$, $c_2\geq (x_1+c_1)-x_2$, $c_3\geq (x_2+c_2)-x_3$,
etc.  Let us denote $y_1=x-c_1$, $y_2=x_1+c_1-c_2$, $y_3=x_2+c_2-c_3,\dots,$
$y_m=x_{m-1}+c_{m-1}-c_m$, $y_{m+1}=x_m+c_m$.  Then all the above
inequalities are equivalent to the interlacing conditions~(\ref{eq:interlace}).
\endproof

Due to Lemma~\ref{le:odot} and Proposition~\ref{pr:e-pieri},
Theorem~\ref{th:main} and Corollary~\ref{cor:rule}
would follow from Theorem~\ref{th:associative},
which says that~$M$ is an associative
operation.  
The proof of associativity given in Section~\ref{sec:associativity}
is based on a Yang-Baxter type relation for the 
scattering matrices $R_{ij}(c)$ (Section~\ref{sec:yb}) and
on the construction of certain polyhedral cones in the space of the 
parameters $c_{ij}$ (Section~\ref{sec:cone}).

\section{Yang-Baxter Equation and Tetrahedron Equation}
\label{sec:yb}
\neweq

As we mentioned before, for distinct $i,j,k,l$, the endomorphism
$R_{ij}(c_{ij})$ commutes with $R_{kl}(c_{kl})$.
Thus $R_\a(C)$ does not change when we apply a $2$-move to the reduced
word~$\a$.  The relations that involve $3$-moves are less trivial.

\begin{theorem} 
\label{th:rrr} 
The operators $R_{12}(c_{12})$, $R_{13}(c_{13})$, and $R_{23}(c_{23})$ 
acting on~$E^{\otimes 3}$ satisfy the relation
\begin{equation} 
R_{23}(c_{23})R_{13}(c_{13})R_{12}(c_{12}) =
R_{12}(c_{12}')R_{13}(c_{13}')R_{23}(c_{23}'),
\label{eq:rrr}
\end{equation} 
where~$c_{12}, c_{13}, c_{23}$ are arbitrary parameters and
$c_{12}', c_{13}', c_{23}'$ are given by
\begin{equation}
\begin{array}{l}
c_{12}'=\min(c_{12},c_{13}-c_{23}),\\[.05in]
c_{13}'=c_{12}+c_{23},\\[.05in]
c_{23}'=\max(c_{23},c_{13}-c_{12}).
\end{array}
\label{eq:transform}
\end{equation}
Moreover, for fixed $c_{12}, c_{13}, c_{23}$ 
the collection $c_{12}', c_{13}', c_{23}'$ defined 
by~{\rm (\ref{eq:transform})} is a unique collection of parameters
such that~{\rm (\ref{eq:rrr})} holds identically.
\end{theorem}

The following two wiring diagrams related by a $3$-move 
illustrate the statement of the theorem.
\psset{curvature=0 0 0}
\psset{xunit=.4cm,yunit=.4cm}
\begin{center}
\pspicture[.1](0,-7)(28,11)

  \pscurve[linewidth=.5pt, showpoints=false]{<-}%
          (2,-1)(2,0)(2,1)(3,2)(4,3)(5,4)(6,5)(6,6)(6,7)(6,8)(6,9)
  \pscurve[linewidth=.5pt, showpoints=false]{<-}%
          (4,-1)(4,0)(4,1)(3,2)(2,3)(2,4)(2,5)(3,6)(4,7)(4,8)(4,9)
  \pscurve[linewidth=.5pt, showpoints=false]{<-}%
          (6,-1)(6,0)(6,1)(6,2)(6,3)(5,4)(4,5)(3,6)(2,7)(2,8)(2,9)
  \qdisk(3,2){1.5pt}
  \qdisk(5,4){1.5pt}
  \qdisk(3,6){1.5pt}
  \rput[u](1.3,6){$c_{12}$}
  \rput[u](6.7,4){$c_{13}$}
  \rput[u](1.3,2){$c_{23}$}
  \rput[u](2,9.7){$x_1$}
  \rput[u](4,9.7){$x_2$}
  \rput[u](6,9.7){$x_3$}
  \rput[u](2,-1.7){$y_3$}
  \rput[u](4,-1.7){$y_2$}
  \rput[u](6,-1.7){$y_1$}

  \psline[linewidth=.2pt]{->}(12,4)(16,4)

  \pscurve[linewidth=.5pt, showpoints=false]{<-}%
          (22,-1)(22,0)(22,1)(22,2)(22,3)(23,4)(24,5)(25,6)(26,7)(26,8)(26,9)
  \pscurve[linewidth=.5pt, showpoints=false]{<-}%
          (24,-1)(24,0)(24,1)(25,2)(26,3)(26,4)(26,5)(25,6)(24,7)(24,8)(24,9)
  \pscurve[linewidth=.5pt, showpoints=false]{<-}%
          (26,-1)(26,0)(26,1)(25,2)(24,3)(23,4)(22,5)(22,6)(22,7)(22,8)(22,9)
  \qdisk(25,2){1.5pt}
  \qdisk(23,4){1.5pt}
  \qdisk(25,6){1.5pt}
  \rput[u](26.7,6){$c_{23}'$}
  \rput[u](21.3,4){$c_{13}'$}
  \rput[u](26.7,2){$c_{12}'$}
  \rput[u](22,9.7){$x_1$}
  \rput[u](24,9.7){$x_2$}
  \rput[u](26,9.7){$x_3$}
  \rput[u](22,-1.7){$y_3$}
  \rput[u](24,-1.7){$y_2$}
  \rput[u](26,-1.7){$y_1$}

  \rput[u](14,-5){$%
        \begin{array}{rcccc}
        y_1&=&x_3+c_{13}&=&x_3+c_{23}'+c_{12}',\\[.05in]
        y_2&=&x_2+c_{12}-c_{13}+c_{23}&=&x_2-c_{23}'+c_{13}'-c_{12}',\\[.05in]
        y_3&=&x_1-c_{12}-c_{23}&=&x_1-c_{13}'
        \end{array}%
$}

\endpspicture
\end{center}

\begin{remark} {\rm
If $c_{13}=c_{12}+c_{23}$ then 
$c'_{12}=c_{12}$, $c'_{13}=c_{13}$, and $c'_{23}=c_{23}$.
In this case the equation~(\ref{eq:rrr}) becomes the famous
{\it quantum Yang-Baxter equation\/} with two parameters, which is 
well-known to the informed reader in the form
$R_{12}(u)R_{13}(u+v)R_{23}(v)=R_{23}(v)R_{13}(u+v)R_{12}(u)$.
}
\end{remark}

\proofof{Theorem~\ref{th:rrr}}
The operator $R_{23}(c_{23})R_{13}(c_{13})R_{12}(c_{12})$ maps 
the basis vector $e_{x_1}\otimes e_{x_2}\otimes e_{x_3}$ 
either to $e_{x_3 + c_{13}}\otimes e_{x_2+c_{12}-c_{13}+c_{23}}\otimes 
e_{x_1-c_{12}-c_{23}}$ if
\begin{equation}
\left\{
\begin{array}{l}
   c_{12} \geq x_1 - x_2, \\[.05in]
   c_{13} \geq (x_2 + c_{12}) - x_3, \\[.05in]
   c_{23} \geq (x_1 - c_{12}) - (x_2 + c_{12} - c_{13}),
\end{array}
\right.
\label{eq:ccc1}
\end{equation}
or to zero otherwise.
Likewise, the operator $R_{12}(c_{12}')R_{13}(c_{13}')R_{23}(c_{23}')$ maps 
$e_{x_1}\otimes e_{x_2}\otimes e_{x_3}$  either to 
$e_{x_3+c_{23}'+c_{12}'}\otimes e_{x_2-c_{23}'+c_{13}'-c_{12}'}\otimes 
e_{x_1-c_{13}'}$ if
\begin{equation}
\left\{
\begin{array}{l}
   c_{23}' \geq x_2 - x_3, \\[.05in]
   c_{13}' \geq x_1 -(x_2 - c_{23}'), \\[.05in]
   c_{12}' \geq (x_2 - c_{23}' + c_{13}') - (x_3 - c_{23}'),
\end{array}
\right.
\label{eq:ccc2}
\end{equation}
or to zero otherwise.
These two operators are equal if and only if
\begin{equation}
\begin{array}{l}
c_{13}=c_{12}'+c_{23}',\\[.05in]
c_{12}-c_{13}+c_{23}=-c_{23}'+c_{13}'-c_{12}',\\[.05in]
c_{12}+c_{23}=c_{13}'.
\end{array}
\label{eq:ccc3}
\end{equation}
(the second identity is the difference of two others)
and for any $x_1,x_2,x_3$ the condition~(\ref{eq:ccc1}) is 
equivalent to the condition~(\ref{eq:ccc2}).
We can write these two sets of inequalities in a more compact form
as follows:
$$
\left\{
\begin{array}{l}
 \min(c_{12},c_{23}+2c_{12}-c_{13})\geq x_1 -x_2,\\[.05in]
 c_{13}-c_{12}\geq x_2-x_3,
\end{array}
\right.
\textrm{is equivalent to}
\left\{
\begin{array}{l}
 c_{13}'-c_{23}'\geq x_1-x_2,\\[.05in]
 \min(c_{23}',c_{12}+2c_{23}'-c_{13}')\geq x_2-x_3.
\end{array}
\right.
$$
Thus 
$$
\begin{array}{l}
 \min(c_{12},c_{23}+2c_{12}-c_{13}) = c_{13}'-c_{23}',\\[.05in]
 c_{13}-c_{12} = \min(c_{23}',c_{12}+2c_{23}'-c_{13}').
\end{array}
$$
These two identities together with~(\ref{eq:ccc3}) are equivalent to
the relations~(\ref{eq:transform}).
\endproof

It follows for Theorem~\ref{th:rrr} that, for $i<j<k$,
the operators $R_{ij}(c_{ij})$, $R_{ik}(c_{ik})$, and $R_{jk}(c_{jk})$ 
acting on~$E^{\otimes n}$ satisfy the relation
$ R_{jk}(c_{jk})R_{ik}(c_{ik})R_{ij}(c_{ij}) =
R_{ij}(c_{ij}')R_{ik}(c_{ik}')R_{jk}(c_{jk}')$,
where
\begin{equation}
\begin{array}{l}
c_{ij}'=\min(c_{ij},c_{ik}-c_{jk}),\\[.05in]
c_{ik}'=c_{ij}+c_{jk},\\[.05in]
c_{jk}'=\max(c_{jk},c_{ik}-c_{ij}).
\end{array}
\label{eq:transform_general}
\end{equation}
The inverse transformation $(c_{ij}',c_{ik}',c_{jk}')\to
(c_{ij},c_{ik},c_{jk})$ is given by similar formulas
\begin{equation}
\begin{array}{l}
c_{ij}=\max(c_{ij}',c_{ik}'-c_{jk}'),\\[.05in]
c_{ik}=c_{ij}'+c_{jk}',\\[.05in]
c_{jk}=\min(c_{jk}',c_{ik}'-c_{ij}').
\end{array}
\label{eq:inverse_transform}
\end{equation}

We will denote by $\Z^{I(w)}$  the set of all collections of integers
parameters $C=(c_{pq})$ with $(p,q)\in I(w)$, $p<q$.  
For $i<j<k$ such that $(i,j), (i,k), (j,k)\in
I(w)$,  we denote by $T_{ijk}$ the local transformation of parameters
$$
\begin{array}{l}
T_{ijk}:\,\Z^{I(w)}\longrightarrow \Z^{I(w)}\\[.15in]
T_{ijk}:\,(c_{pq})\longmapsto (c_{pq}')
\end{array}
$$
where the $c_{pq}'$ are given by formulas~(\ref{eq:transform_general}) 
for $p,q\in\{i,j,k\}$ and $c_{pq}'=c_{pq}$ otherwise.

For any two reduced words~$\a,\b\in\R(w)$ of a permutation $w\in S_n$,  we
will define a transition map $T_\a^\b:\Z^{I(w)}\to\Z^{I(w)}$ as a
composition of local transformation maps~$T_{ijk}$.  If
$\a=(\dots,a,b,\dots)$ and $\a'=(\dots,b,a,\dots)$, $|a-b|\geq 2$, are two
reduced words for~$w$ related by a $2$-move, then $T_\a^{\a'}$ is the identity
map.   If $\a=(\dots,a, a+1, a,\cdots)$ and $\a'=(\dots,a+1,a,a+1,\cdots)$ are
two reduced words related by a $3$-move, then  the corresponding reflection
orderings of~$I(w)$ differ only in three places:
$\cdots< (j,k)< (i,k)< (i,j)< \cdots$ and 
$\cdots< (i,j)< (i,k)< (j,k)< \cdots$ for certain
$i<j<k$.  In this case we define $T_\a^{\a'}=T_{ijk}$ and $T_{\a'}^\a=T_{ijk}^{-1}$.
In general, we choose a chain of reduced words 
$\a,\a^1,\a^2,\dots,\a^k,\b\in\R(w)$ 
that interpolates between $\a$ and $\b$ such that any two adjacent
words are related by a $2$-~or $3$-move.  
Then we define $T_\a^\b= T_{\a^k}^{\b}\cdots T_{\a^1}^{\a^2}\,T_\a^{\a^1}$.

It follows from the uniqueness part of Theorem~\ref{th:rrr} that the transition
map $T_\a^{\b}$ does not depend upon a choice of path of $2$- and $3$-moves
joining the reduced words $\a$ and $\b$.  Let us remark that
this property amounts to
verifying that the local transformation maps $T_{ijk}$ satisfy the following
tetrahedron equation.

\begin{theorem}  {\rm (Tetrahedron equation)} \
Let $\wnot$ be the longest element in~$S_4$.  The following identity for the 
compositions of maps $\Z^{I(\wnot)}\to\Z^{I(\wnot)}$ holds:
$$ 
T_{123} T_{124} T_{134} T_{234} = T_{234} T_{134} T_{124} T_{123}.
$$
\end{theorem}

It is left as an exercise for the reader to verify directly that the
local transformation maps $T_{ijk}$ satisfy the tetrahedron equation.

\medskip
Recall that $R_\a(C)$ is the composition of scattering matrices
defined by~(\ref{eq:rac}).  It is immediately clear from Theorem~\ref{th:rrr} that
$ R_\a(C)=R_{\a'}(T_\a^{\a'}(C))$ if $\a$ and $\a'$ are related by a $2$~or
$3$-move.  Thus, in general, we have
\begin{equation}
R_\a(C)=R_{\b}(T_\a^{\b}(C))
\label{eq:Rab}
\end{equation}
for any two reduced words $\a$ and $\b$ for~$w$ and any collection of parameters
 $C\in\Z^{I(w)}$.

\section{Principal Cones}                    
\label{sec:cone}
\neweq

Let $\a$ be a reduced word of a permutation~$w\in S_n$.
In this section we construct and study a certain polyhedral cone $\C_\a$ in 
the space $\Z^{I(w)}$.  In the case when $w=\wnot$ is the longest
permutation in $S_n$, the cone $\C_\a$ is exactly the cone
of Kashiwara's parametrizations of dual canonical bases for 
$U_q(sl_n)$.  It is the {\it string cone} in the terminology of 
Berenstein and Zelevinsky~\cite{BZ0}.
The explicit description of $\C_\a$ gives an answer
to a question posed in~\cite{BZ0}.

\subsection{Rigorous paths and statements of results}

Let us fix a reduced word $\a\in\R(w)$ and an integer $0\leq s \leq n$.  
We construct an oriented graph
$G(\a,s)$ from the wiring diagram corresponding to $\a$ as follows.  Denote by 
$v_{ij}$ the vertex of the wiring diagram which is the intersection of the
$i$-th and $j$-th pseudo-lines.  The vertex set of the graph $G(\a,s)$ is
composed of  the vertices $v_{ij}$ together with $2n$ boundary vertices:
$U_1,\dots,U_n$ that mark the upper ends of pseudo-lines from left to right and
$L_1,\dots,L_n$, that mark the lower ends of pseudo-lines from left to right.
Notice that $U_i$ is the upper end of the $i$-th pseudo-line.  We orient
downward the $s$ pseudo-lines of the wiring diagram whose lower ends are
labelled $L_1,\dots,L_s$ and we orient upward the remaining $n-s$ pseudo-lines
whose lower ends are labelled $L_{s+1},\dots,L_n$.  Two vertices are connected
by an edge in the graph $G(\a,s)$ if they are adjacent vertices on the same
pseudo-line.  Directions of edges in $G(\a,s)$ agree with directions of the
corresponding pseudo-lines.  For example, the graph $G(121, 2)$ is shown on the
picture below.

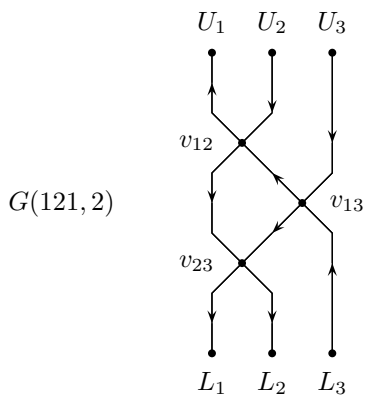
\psset{curvature=0 0 0}
\psset{xunit=.4cm,yunit=.4cm}
\begin{center}
\pspicture[.1](0,-3)(7,11)

  \psline[linewidth=.5pt, showpoints=false]{-}%
          (2,-1)(2,0)(2,1)(3,2)(4,3)(5,4)(6,5)(6,6)(6,7)(6,8)(6,9)
  \psline[linewidth=.5pt, showpoints=false]{-}%
          (4,-1)(4,0)(4,1)(3,2)(2,3)(2,4)(2,5)(3,6)(4,7)(4,8)(4,9)
  \psline[linewidth=.5pt, showpoints=false]{-}%
          (6,-1)(6,0)(6,1)(6,2)(6,3)(5,4)(4,5)(3,6)(2,7)(2,8)(2,9)

  \psline[linewidth=.5pt, showpoints=false]{-}(2,-1)(2,0)
  \psline[linewidth=.5pt, showpoints=false]{<-}(2,0)(2,1)(4,3)
  \psline[linewidth=.5pt, showpoints=false]{<-}(4,3)(6,5)(6,6)
  \psline[linewidth=.5pt, showpoints=false]{<-}(6,6)(6,9)

  \psline[linewidth=.5pt, showpoints=false]{-}(4,-1)(4,0)
  \psline[linewidth=.5pt, showpoints=false]{<-}(4,0)(4,1)(2,3)(2,4)
  \psline[linewidth=.5pt, showpoints=false]{<-}(2,4)(2,5)(4,7)
  \psline[linewidth=.5pt, showpoints=false]{<-}(4,7)(4,9)

  \psline[linewidth=.5pt, showpoints=false]{->}(6,-1)(6,2)
  \psline[linewidth=.5pt, showpoints=false]{->}(6,2)(6,3)(4,5)
  \psline[linewidth=.5pt, showpoints=false]{->}(4,5)(2,7)(2,8)
  \psline[linewidth=.5pt, showpoints=false]{-}(2,8)(2,9)

  \qdisk(3,2){1.5pt}
  \qdisk(5,4){1.5pt}
  \qdisk(3,6){1.5pt}
  \qdisk(2,-1){1.5pt}
  \qdisk(4,-1){1.5pt}
  \qdisk(6,-1){1.5pt}
  \qdisk(6,9){1.5pt}
  \qdisk(4,9){1.5pt}
  \qdisk(2,9){1.5pt}
  \rput[u](1.5,6){$v_{12}$}
  \rput[u](6.5,4){$v_{13}$}
  \rput[u](1.5,2){$v_{23}$}
  \rput[u](2,10){$U_1$}
  \rput[u](4,10){$U_2$}
  \rput[u](6,10){$U_3$}
  \rput[u](6,-2){$L_3$}
  \rput[u](4,-2){$L_2$}
  \rput[u](2,-2){$L_1$}

  \rput[u](-3,4) {$G(121,2)$}
\endpspicture
\end{center}

An oriented path in the graph $G(\a,s)$ is a sequence of vertices 
$v_0,\dots,v_l$ connected by the oriented edges $v_0\to v_1,\ v_1\to v_2,\ 
\dots,\ v_{l-1}\to v_l$.  Notice that the graph $G(\a,s)$ is acyclic,
i.e., there is no closed oriented cycle in the graph.
Thus there are finitely many oriented paths in $G(\a,s)$.
We say that an oriented path $v_0\to v_1\to  \dots \to v_l$ 
is {\it rigorous} if it satisfies the following condition:
There are no three adjacent vertices $v_{a} \to v_{a+1} \to v_{a+2}$ 
in the path
such that $v_a,v_{a+1},v_{a+2}$ belong to the same $i$-th pseudo-line,
$v_{a+1}$ is the intersection of the $i$-th and $j$-th pseudo-line, and either
$i<j$ and both $i$-th and $j$-th pseudo-lines
are oriented upwards, or $i>j$ and the $i$-th and $j$-th pseudo-lines are 
oriented downwards.  
In other words, a path is rigorous if and only if it avoids the following two
fragments:

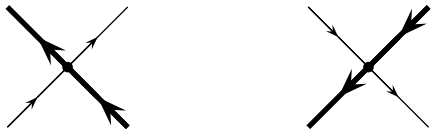
\psset{curvature=0 0 0}
\psset{xunit=.4cm,yunit=.4cm}
\begin{center}
\pspicture[.1](0,0)(14,4)

  \psline[linewidth=.5pt, showpoints=false]{->}(0,0)(1,1)
  \psline[linewidth=.5pt, showpoints=false]{->}(1,1)(3,3)
  \psline[linewidth=.5pt, showpoints=false]{-}(0,0)(4,4)
  \psline[linewidth=2pt, showpoints=false]{->}(4,0)(3,1)
  \psline[linewidth=2pt, showpoints=false]{->}(3,1)(1,3)
  \psline[linewidth=2pt, showpoints=false]{-}(4,0)(0,4)
  \qdisk(2,2){2pt}

  \psline[linewidth=2pt, showpoints=false]{-}(10,0)(14,4)
  \psline[linewidth=2pt, showpoints=false]{<-}(11,1)(13,3)
  \psline[linewidth=2pt, showpoints=false]{<-}(13,3)(14,4)
  \psline[linewidth=.5pt, showpoints=false]{-}(14,0)(10,4)
  \psline[linewidth=.5pt, showpoints=false]{<-}(13,1)(11,3)
  \psline[linewidth=.5pt, showpoints=false]{<-}(11,3)(10,4)
  \qdisk(12,2){2pt}

\endpspicture
\end{center}
Here the thick lines show path fragments and the thin lines show the
pseudo-lines they intersect.

For example, in the graph $G(121,2)$ shown above 
all paths connecting boundary vertices are rigorous except the following
two paths: $L_3\to v_{13}\to v_{23}\to L_1$ and 
$U_3\to v_{13}\to v_{23} \to L_1$.

Let $P=(v_0\to v_1\to  \cdots \to v_l)$ be a rigorous path 
connecting two boundary vertices $v_0$ and~$v_l$.  Suppose that the edge
$v_{r-1}\to v_{r}$ is on the $i_r$-th pseudo-line, for $r=1,\dots,l$.
We denote by $c_P$ the expression
\begin{equation}
c_P=c_{i_1 i_2}+c_{i_2 i_3}+\cdots+c_{i_{l-1}\,i_l},
\label{eq:c_P}
\end{equation}
where we assume that $c_{ii}=0$ and for $i>j$ the coefficients $c_{ij}$ are given by
$c_{ij}=-c_{ji}$.

\begin{definition} {\rm
For a reduced word $\a\in\R(w)$, we define the {\it principal cone}
$\C_\a$ as the the polyhedral cone in the the integer lattice $\Z^{I(w)}$
of collections $C=(c_{ij})$ given by the inequalities 
$c_P\geq 0$ for all rigorous paths $P$ in the graph $G(\a,s)$ from the vertex
$L_{s+1}$ to $L_s$, for $1\leq s\leq n-1$.}
\label{def:princ}
\end{definition}

\begin{theorem}
For any two reduced words $\a,\, \b\in\R(w)$, the transition map $T_\a^\b$
bijectively maps the cone $\C_\a$ to the cone $\C_\b$.
\label{th:cones-bijectively}
\end{theorem}

\begin{example} {\rm
To illustrate the definition and the theorem we describe the principal cones 
for two reduced decompositions of~$\wnot$ in~$S_3$.
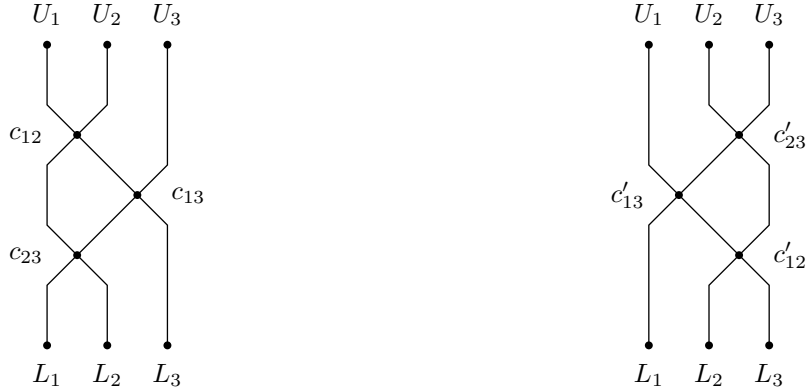
\psset{curvature=0 0 0}
\psset{xunit=.4cm,yunit=.4cm}
\begin{center}
\pspicture[.1](0,-5)(28,11)

  \psline[linewidth=.5pt, showpoints=false]{-}%
          (2,-1)(2,0)(2,1)(3,2)(4,3)(5,4)(6,5)(6,6)(6,7)(6,8)(6,9)
  \psline[linewidth=.5pt, showpoints=false]{-}%
          (4,-1)(4,0)(4,1)(3,2)(2,3)(2,4)(2,5)(3,6)(4,7)(4,8)(4,9)
  \psline[linewidth=.5pt, showpoints=false]{-}%
          (6,-1)(6,0)(6,1)(6,2)(6,3)(5,4)(4,5)(3,6)(2,7)(2,8)(2,9)
  \qdisk(3,2){1.5pt}
  \qdisk(5,4){1.5pt}
  \qdisk(3,6){1.5pt}
  \rput[u](1.3,6){$c_{12}$}
  \rput[u](6.7,4){$c_{13}$}
  \rput[u](1.3,2){$c_{23}$}

  \rput[u](2,10){$U_1$}
  \rput[u](4,10){$U_2$}
  \rput[u](6,10){$U_3$}

  \qdisk(2,9){1.5pt}
  \qdisk(4,9){1.5pt}
  \qdisk(6,9){1.5pt}

  \qdisk(2,-1){1.5pt}
  \qdisk(4,-1){1.5pt}
  \qdisk(6,-1){1.5pt}

  \rput[u](2,-2){$L_1$}
  \rput[u](4,-2){$L_2$}
  \rput[u](6,-2){$L_3$}
  
  \psline[linewidth=.5pt, showpoints=false]{-}%
          (22,-1)(22,0)(22,1)(22,2)(22,3)(23,4)(24,5)(25,6)(26,7)(26,8)(26,9)
  \psline[linewidth=.5pt, showpoints=false]{-}%
          (24,-1)(24,0)(24,1)(25,2)(26,3)(26,4)(26,5)(25,6)(24,7)(24,8)(24,9)
  \psline[linewidth=.5pt, showpoints=false]{-}%
          (26,-1)(26,0)(26,1)(25,2)(24,3)(23,4)(22,5)(22,6)(22,7)(22,8)(22,9)
  \qdisk(25,2){1.5pt}
  \qdisk(23,4){1.5pt}
  \qdisk(25,6){1.5pt}
  \rput[u](26.7,6){$c_{23}'$}
  \rput[u](21.3,4){$c_{13}'$}
  \rput[u](26.7,2){$c_{12}'$}

  \rput[u](22,10){$U_1$}
  \rput[u](24,10){$U_2$}
  \rput[u](26,10){$U_3$}

  \qdisk(22,9){1.5pt}
  \qdisk(24,9){1.5pt}
  \qdisk(26,9){1.5pt}

  \qdisk(22,-1){1.5pt}
  \qdisk(24,-1){1.5pt}
  \qdisk(26,-1){1.5pt}

  \rput[u](22,-2){$L_1$}
  \rput[u](24,-2){$L_2$}
  \rput[u](26,-2){$L_3$}

  \rput[u](3,-4) {$\C_{121}=\{c_{12} \geq 0, \ c_{13} \geq c_{23} \geq 0\}$}
  \rput[u](23,-4){$\C_{212}=\{c_{23}'\geq 0, \ c_{13}'\geq c_{12}'\geq 0\}$}

\endpspicture
\end{center}
Indeed, for $\a=121$ we have the rigorous paths the $L_2\to v_{23}\to L_1$,
$L_3\to v_{13}\to v_{23}\to L_2$, and $L_3\to v_{13}\to v_{12}\to v_{23}\to L_2$.
Analogously, for $\a=212$ we have the rigorous paths $L_2\to v_{12}\to v_{23}\to v_{13}
\to L_1$, $L_2\to v_{12}\to v_{13}\to L_1$, and $L_3\to v_{12}\to L_2$.
One can easily verify that the transformation map $T_{123}$ maps the 
cone $\C_{121}$ into the cone $\C_{212}$.
}
\end{example}

In the case when $\a\in\R(\wnot)$ is a reduced word for the longest permutation
in $S_n$ there are two alternative descriptions of the principal cone $\C_\a$.

\begin{theorem}
For a reduced word $\a\in\R(\wnot)$, the principal cone
$\C_\a$ is the
set of all collections $C=(c_{ij})\in \Z^{I(\wnot)}$ such that for any reduced
word $\b\in\R(\wnot)$ all entries~$c_{ij}'$ of the collection 
$C'=(c_{ij}')=T_\a^{\b}(C)\in \Z^{I(w)}$ are nonnegative.
\label{th:princ-w0}
\end{theorem}

For a reduced word $\b\in \R(\wnot)$, let $low(\b)$ denote the pair $(i,j)$,
$1\leq i<j\leq n$, such that the lowest node of the wiring diagram of $\b$ 
is the crossing of $i$-th and $j$-th pseudo-lines.  
(It is clear that $j=i+1$.) 
For example, $low(121)=(2,3)$ and $low(212)=(1,2)$.

The principal cone can be described by a weaker set of conditions
as follows.

\begin{theorem}
\label{th:low}
For a reduced word $\a\in\R(\wnot)$, the principal cone $\C_\a$ is the set 
of all collections $C=(c_{ij})\in\Z^{I(\wnot)}$
such that for any reduced word $\b\in\R(\wnot)$ the lowest entry
$c'_{low(\b)}$ of $C'=(c_{ij}')=T_\a^\b(C)$ is nonnegative.
\end{theorem}

\begin{remark} {\rm
In the case of $\wnot$, we can choose either of the descriptions from
Theorems~\ref{th:princ-w0} or~\ref{th:low} as the definition of 
the principal cone.  Then Theorem~\ref{th:cones-bijectively} would
become trivial.  But this would obscure the fact that the principal 
cone is actually a polyhedral cone.  }
\end{remark}

Before we proceed, let us consider several examples of principal cones.

\begin{example} {\rm
Let $\anot=(1,2,1,3,2,1,\dots,n-1,n-2,\dots,1)\in \R(\wnot)$ 
be the lexicographically minimal reduced word for the longest permutation. 
By Definition~\ref{def:princ}, the principal cone $\C_\anot$ is given 
by the inequalities: 
\begin{equation} 
c_{12}\geq 0\,, \quad  c_{13}\geq c_{23}\geq
0\,, \quad c_{14}\geq c_{24}\geq c_{34} \geq 0\,, \quad c_{15}\geq c_{25}\geq
c_{35} \geq c_{45} \geq 0\,, \ \dots 
\label{eq:string-cc} 
\end{equation}
Indeed, in this case all inequalities $c_P\geq 0$ in Definition~\ref{def:princ}
are of the form $c_{s\,s+1}\geq 0$ and $c_{si}-c_{s+1\,i}\geq 0$ for $i > s+1$.  
}
\label{ex:anot}
\end{example}

Berenstein and Zelevinsky~\cite{BZ0} studied the string cone of Kashiwara's
parametrizations of dual canonical basis for $U_q(sl_n)$.  This is a cone
$\tilde\C_\a$ in the ${n\choose 2}$-dimensional space of strings $C=(c_{ij})$
that depends upon a choice of reduced word $\a\in\R(\wnot)$ for the longest
permutation.  It follows from the definitions that if $\a$ and $\a'$ differ by a
$2$- or $3$-move then $\tilde\C_{\a'} = T_\a^{\a'}(\tilde\C_\a)$.  Thus string
cones $\tilde\C_\a$ transform according to the transition maps $T_\a^\b$.  The
string cone was explicitly calculated in~\cite{BZ0} for the lexicographically
minimal reduced word $\anot$.  In this case $\tilde\C_\anot$ is given by the
inequalities~(\ref{eq:string-cc}).  Theorem~\ref{th:cones-bijectively} and
Example~\ref{ex:anot} imply the following statement.

\begin{corollary}
{\rm (String cones)} \
For a reduced word $\a\in\R(\wnot)$, the principal cone $\C_\a$
is exactly the string cone~$\tilde\C_\a$.
\label{cor:string}
\end{corollary}

Definition~\ref{def:princ} gives an explicit description of the string
cone~$\tilde\C_\a$.  This settles the problem of describing the string cones 
for any reduced word~$\a\in\R(\wnot)$.

\begin{example} {\rm
This example is related to our construction of the $*$-product in 
Section~\ref{sec:scattering}.
Recall that the permutation $w(m,n): i \mapsto i+n \pmod{m+n}$ in $S_{m+n}$
has a unique reduced decomposition up to $2$-moves.  For example, for
$m=4$ and $n=3$ we have the following wiring diagram:

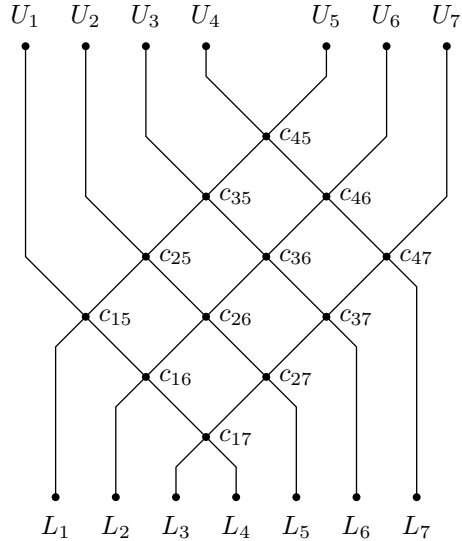
\psset{curvature=0 0 0}
\psset{xunit=.4cm,yunit=.4cm}
\begin{center}
\pspicture[.1](0,-1)(16,18)

  \pscurve[linewidth=.5pt, showpoints=false]{-}%
          (2,1)(2,6)(11,15)(11,16)
  \pscurve[linewidth=.5pt, showpoints=false]{-}%
          (4,1)(4,4)(13,13)(13,16)
  \pscurve[linewidth=.5pt, showpoints=false]{-}%
          (6,1)(6,2)(15,11)(15,16)

  \pscurve[linewidth=.5pt, showpoints=false]{-}%
          (8,1)(8,2)(1,9)(1,16)
  \pscurve[linewidth=.5pt, showpoints=false]{-}%
          (10,1)(10,4)(3,11)(3,16)
  \pscurve[linewidth=.5pt, showpoints=false]{-}%
          (12,1)(12,6)(5,13)(5,16)
  \pscurve[linewidth=.5pt, showpoints=false]{-}%
          (14,1)(14,8)(7,15)(7,16)

  \qdisk(7,3){1.5pt}
  \qdisk(9,5){1.5pt}
  \qdisk(11,7){1.5pt}
  \qdisk(13,9){1.5pt}
  \qdisk(5,5){1.5pt}
  \qdisk(7,7){1.5pt}
  \qdisk(9,9){1.5pt}
  \qdisk(11,11){1.5pt}
  \qdisk(3,7){1.5pt}
  \qdisk(5,9){1.5pt}
  \qdisk(7,11){1.5pt}
  \qdisk(9,13){1.5pt}

  \rput[u](8,3){$c_{17}$}
  \rput[u](6,5){$c_{16}$}
  \rput[u](4,7){$c_{15}$}
  \rput[u](10,5){$c_{27}$}
  \rput[u](8,7){$c_{26}$}
  \rput[u](6,9){$c_{25}$}
  \rput[u](12,7){$c_{37}$}
  \rput[u](10,9){$c_{36}$}
  \rput[u](8,11){$c_{35}$}
  \rput[u](14,9){$c_{47}$}
  \rput[u](12,11){$c_{46}$}
  \rput[u](10,13){$c_{45}$}

  \qdisk(2,1){1.5pt}
  \qdisk(4,1){1.5pt}
  \qdisk(6,1){1.5pt}
  \qdisk(8,1){1.5pt}
  \qdisk(10,1){1.5pt}
  \qdisk(12,1){1.5pt}
  \qdisk(14,1){1.5pt}

  \qdisk(1,16){1.5pt}
  \qdisk(3,16){1.5pt}
  \qdisk(5,16){1.5pt}
  \qdisk(7,16){1.5pt}
  \qdisk(11,16){1.5pt}
  \qdisk(13,16){1.5pt}
  \qdisk(15,16){1.5pt}

  \rput[u](1,17){$U_1$}
  \rput[u](3,17){$U_2$}
  \rput[u](5,17){$U_3$}
  \rput[u](7,17){$U_4$}

  \rput[u](11,17){$U_5$}
  \rput[u](13,17){$U_6$}
  \rput[u](15,17){$U_7$}

  \rput[u](2,0){$L_1$}
  \rput[u](4,0){$L_2$}
  \rput[u](6,0){$L_3$}
  \rput[u](8,0){$L_4$}
  \rput[u](10,0){$L_5$}
  \rput[u](12,0){$L_6$}
  \rput[u](14,0){$L_7$}
\endpspicture
\end{center}
By Definition~\ref{def:princ}, the corresponding principal cone $\C_{(m,n)}=\C_\a$ 
is given by the inequalities $-c_{ik} + c_{ij}\geq 0$ for $i\leq m < j<k$;
$c_{1\,m+n}\geq 0$; and $c_{jk}-c_{ik}\geq 0$ for $i<j\leq m<k$.
These are exactly the conditions~(\ref{eq:conditions}) on the parameters
in the sum~(\ref{eq:mmn}).  Thus the $*$-product in $T(E)$ can be written
as the sum 
$$
e_{x_1\dots x_m} * e_{y_1\dots y_n} =  \sum_{C\in\C_{(m,n)}}R_{(m,n)}(C)
\cdot (e_{x_1\dots x_m} \otimes e_{y_1\dots y_n}).
$$
} 
\label{ex:rmn}
\end{example}

\begin{example} {\rm 
Let us also illustrate the definitions by the following example for 
the reduced decomposition $s_2\, s_1\, s_2\, s_3\, s_2\, s_1$ of 
the longest element $\wnot\in S_4$.

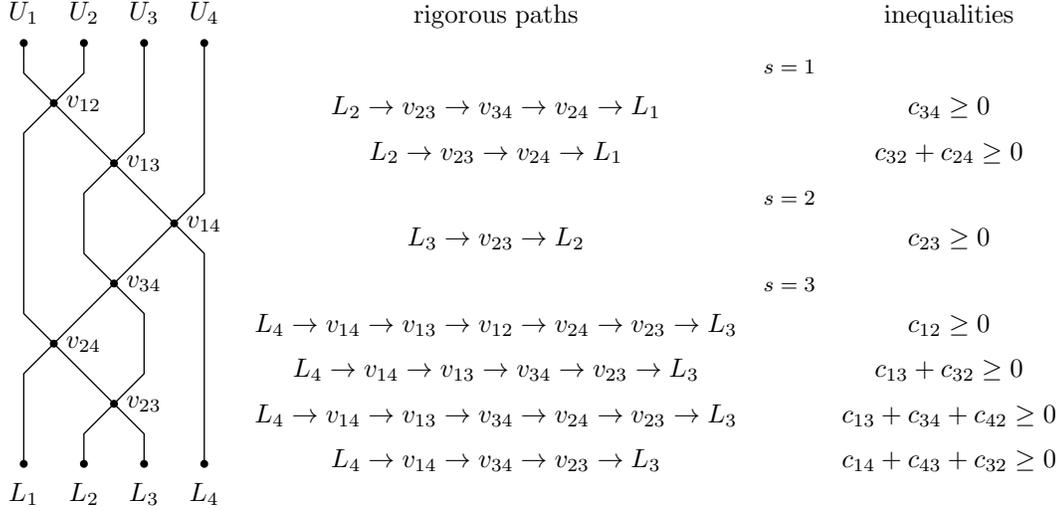
\psset{curvature=0 0 0}
\psset{xunit=.4cm,yunit=.4cm}
\begin{center}
\pspicture[.1](0,-0.5)(37,16.5)

  \psline[linewidth=.5pt, showpoints=false]{-}%
          (2,1)(2,2) (2,3)(2,4) (3,5)(4,6) (5,7)(6,8) (7,9)(8,10) 
          (8,11)(8,12) (8,13)(8,14) (8,15)
  \psline[linewidth=.5pt, showpoints=false]{-}%
          (4,1)(4,2) (5,3)(6,4) (6,5)(6,6) (5,7)(4,8) (4,9)(4,10) 
          (5,11)(6,12) (6,13)(6,14) (6,15)
  \psline[linewidth=.5pt, showpoints=false]{-}%
          (6,1)(6,2) (5,3)(4,4) (3,5)(2,6) (2,7)(2,8) (2,9)(2,10) 
          (2,11)(2,12) (3,13)(4,14) (4,15)
  \psline[linewidth=.5pt, showpoints=false]{-}%
          (8,1)(8,2) (8,3)(8,4) (8,5)(8,6) (8,7)(8,8) (7,9)(6,10) 
          (5,11)(4,12) (3,13)(2,14) (2,15)
  \qdisk(5,3){1.5pt}
  \qdisk(3,5){1.5pt}
  \qdisk(5,7){1.5pt}
  \qdisk(7,9){1.5pt}
  \qdisk(5,11){1.5pt}
  \qdisk(3,13){1.5pt}
  
  \qdisk(2,1){1.5pt}
  \qdisk(4,1){1.5pt}
  \qdisk(6,1){1.5pt}
  \qdisk(8,1){1.5pt}
  \qdisk(2,15){1.5pt}
  \qdisk(4,15){1.5pt}
  \qdisk(6,15){1.5pt}
  \qdisk(8,15){1.5pt}

  \rput[u](2,16){$U_1$}
  \rput[u](4,16){$U_2$}
  \rput[u](6,16){$U_3$}
  \rput[u](8,16){$U_4$}
  \rput[u](2,0){$L_1$}
  \rput[u](4,0){$L_2$}
  \rput[u](6,0){$L_3$}
  \rput[u](8,0){$L_4$}
  \rput[u](4,13){$v_{12}$}
  \rput[u](6,11){$v_{13}$}
  \rput[u](8,9){$v_{14}$}
  \rput[u](6,7){$v_{34}$}
  \rput[u](4,5){$v_{24}$}
  \rput[u](6,3){$v_{23}$}

  \rput[u](23,8.5){$%
\begin{array}{ccc}
\textrm{rigorous paths} && 
\textrm{inequalities}\\[.1in]
&\textrm{\footnotesize $s=1$} &   \\[.05in]
L_2\to v_{23}\to v_{34}\to v_{24}\to L_1   & & 
c_{34}\geq 0\\[.07in]
L_2\to v_{23}\to v_{24}\to L_1   & & 
c_{32}+c_{24}\geq 0 \\[.07in]
&\textrm{\footnotesize $s=2$} &  \\[.05in]
L_3\to v_{23}\to L_2     & & 
c_{23}\geq 0\\[.07in]
&\textrm{\footnotesize $s=3$} &   \\[.05in]
L_4\to v_{14}\to v_{13}\to v_{12}\to v_{24}\to v_{23}\to L_3     & & 
c_{12}\geq 0\\[.07in]
L_4\to v_{14}\to v_{13}\to v_{34}\to v_{23}\to L_3   & & 
c_{13}+c_{32}\geq 0\\[.07in]
L_4\to v_{14}\to v_{13}\to v_{34}\to v_{24}\to v_{23}\to L_3 & & 
c_{13}+c_{34}+c_{42}\geq 0\\[.07in]
L_4\to v_{14}\to v_{34}\to v_{23}\to L_3 & & 
c_{14}+c_{43}+c_{32}\geq 0
\end{array}$}

\endpspicture
\end{center}

In more conventional notation the inequalities defining the 
cone~$\C_{212321}$ can be written as 
$$
\C_{212321}\quad = \quad
\left\{\begin{array}{c}
  c_{12}\geq 0\\[.1in]
  c_{13}\geq c_{23}\geq 0\\[.1in]
  c_{13}+c_{34}\geq c_{24}\geq c_{23}\\[.1in]
  c_{14}\geq c_{23}+c_{34}\\[.1in]
  c_{34}\geq 0
\end{array}
\right\}.
$$
}
\end{example}

\medskip
To prove Theorems~\ref{th:cones-bijectively}, 
\ref{th:princ-w0}, and~\ref{th:low},
we need some extra notation. 
For two boundary vertices
$B$ and $E$ in $\{U_1,\dots,U_n,L_1,\dots,L_n\}$, let 
$$
M_{B,E}^{\a,s} = M_{B,E}^{\a,s}(C) = \min c_P
$$
be the minimum of
expressions~(\ref{eq:c_P}) over all rigorous paths $P$ in the graph $G(\a,s)$
from the vertex $B$ to the vertex $E$, here $C=(c_{ij})$.
(Note that there are finitely many such paths.)
If there are no rigorous paths from $B$ to $E$ in $G(\a,s)$ then
we set $M_{B,E}^{\a,s}=+\infty$.

Using this notation, Definition~\ref{def:princ}
of the principal cone can be written as 
\begin{equation}
\C_\a =\left\{C\in \Z^{I(w)}\mid M_{L_{s+1}, L_s}^{\a,s}(C) \geq 0,\textrm{ for }
s=1,\dots,n-1\right\}.
\label{eq:Ca}
\end{equation}

Theorem~\ref{th:cones-bijectively} is an immediate corollary of the following
more general statement.

\begin{theorem}
\label{th:aa'}
For any two reduced words $\a,\b\in \R(\wnot)$, an integer 
$0\leq s \leq n$, and two boundary vertices $B$ and $E$,
we have $M_{B,E}^{\a,s} (C)=M_{B,E}^{\b,s}(C')$, where $C'=T_\a^\b(C)$.
In other words, the expressions $M_{B,E}^{\a,s}(C)$ are invariant under
the transition maps $T_\a^\b$.
\end{theorem}

\subsection{Proofs}

\proofof{Theorem~\ref{th:aa'}} 
Let us fist verify the statement of the theorem for 
the symmetric group~$S_3$.
In this case we have only two reduced words $121$ and $212$ for $\wnot$.
There are four possible cases $s=0$, $s=1$, $s=2$, and $s=3$.

Let us start with the case $s=3$ when all pseudo-lines are oriented
downwards.  The graphs $G(121,0)$ and $G(212,0)$ are given on 
the picture below.  We also give the transition maps $T_{121}^{212}$
and its inverse $T_{212}^{121}$ for a quick reference, 
cf.~(\ref{eq:transform}) and~(\ref{eq:inverse_transform}).

\psset{curvature=0 0 0}
\psset{xunit=.4cm,yunit=.4cm}
\begin{center}
\pspicture[.1](0,-3)(28,11)

  \psline[linewidth=.5pt, showpoints=false]{-}%
          (2,-1)(2,0)(2,1)(3,2)(4,3)(5,4)(6,5)(6,6)(6,7)(6,8)(6,9)
  \psline[linewidth=.5pt, showpoints=false]{-}%
          (4,-1)(4,0)(4,1)(3,2)(2,3)(2,4)(2,5)(3,6)(4,7)(4,8)(4,9)
  \psline[linewidth=.5pt, showpoints=false]{-}%
          (6,-1)(6,0)(6,1)(6,2)(6,3)(5,4)(4,5)(3,6)(2,7)(2,8)(2,9)

  \psline[linewidth=.5pt, showpoints=false]{-}(2,-1)(2,0)
  \psline[linewidth=.5pt, showpoints=false]{<-}(2,0)(2,1)(4,3)
  \psline[linewidth=.5pt, showpoints=false]{<-}(4,3)(6,5)(6,6)
  \psline[linewidth=.5pt, showpoints=false]{<-}(6,6)(6,9)

  \psline[linewidth=.5pt, showpoints=false]{-}(4,-1)(4,0)
  \psline[linewidth=.5pt, showpoints=false]{<-}(4,0)(4,1)(2,3)(2,4)
  \psline[linewidth=.5pt, showpoints=false]{<-}(2,4)(2,5)(4,7)
  \psline[linewidth=.5pt, showpoints=false]{<-}(4,7)(4,9)

  \psline[linewidth=.5pt, showpoints=false]{-}(6,-1)(6,2)
  \psline[linewidth=.5pt, showpoints=false]{<-}(6,2)(6,3)(4,5)
  \psline[linewidth=.5pt, showpoints=false]{<-}(4,5)(2,7)(2,8)
  \psline[linewidth=.5pt, showpoints=false]{<-}(2,8)(2,9)

  \qdisk(3,2){1.5pt}
  \qdisk(5,4){1.5pt}
  \qdisk(3,6){1.5pt}
  \qdisk(2,-1){1.5pt}
  \qdisk(4,-1){1.5pt}
  \qdisk(6,-1){1.5pt}
  \qdisk(6,9){1.5pt}
  \qdisk(4,9){1.5pt}
  \qdisk(2,9){1.5pt}
  \rput[u](1.3,6){$c_{12}$}
  \rput[u](6.7,4){$c_{13}$}
  \rput[u](1.3,2){$c_{23}$}
  \rput[u](2,10){$U_1$}
  \rput[u](4,10){$U_2$}
  \rput[u](6,10){$U_3$}
  \rput[u](6,-2){$L_3$}
  \rput[u](4,-2){$L_2$}
  \rput[u](2,-2){$L_1$}

  \psline[linewidth=.5pt, showpoints=false]{-}%
          (22,-1)(22,0)(22,1)(22,2)(22,3)(23,4)(24,5)(25,6)(26,7)(26,8)(26,9)
  \psline[linewidth=.5pt, showpoints=false]{-}%
          (24,-1)(24,0)(24,1)(25,2)(26,3)(26,4)(26,5)(25,6)(24,7)(24,8)(24,9)
  \psline[linewidth=.5pt, showpoints=false]{-}%
          (26,-1)(26,0)(26,1)(25,2)(24,3)(23,4)(22,5)(22,6)(22,7)(22,8)(22,9)

  \psline[linewidth=.5pt, showpoints=false]{<-}(22,2)(22,3)
  \psline[linewidth=.5pt, showpoints=false]{<-}(24,5)(25,6)
  \psline[linewidth=.5pt, showpoints=false]{<-}(26,8)(26,9)

  \psline[linewidth=.5pt, showpoints=false]{<-}(24,0)(24,1)
  \psline[linewidth=.5pt, showpoints=false]{<-}(26,4)(26,5)
  \psline[linewidth=.5pt, showpoints=false]{<-}(24,7)(24,8)

  \psline[linewidth=.5pt, showpoints=false]{<-}(26,0)(26,1)
  \psline[linewidth=.5pt, showpoints=false]{<-}(24,3)(23,4)
  \psline[linewidth=.5pt, showpoints=false]{<-}(22,6)(22,7)
  \qdisk(25,2){1.5pt}
  \qdisk(23,4){1.5pt}
  \qdisk(25,6){1.5pt}
  \rput[u](26.7,6){$c_{23}'$}
  \rput[u](21.3,4){$c_{13}'$}
  \rput[u](26.7,2){$c_{12}'$}
  \rput[u](22,10){$U_1$}
  \rput[u](24,10){$U_2$}
  \rput[u](26,10){$U_3$}
  \rput[u](22,-2){$L_1$}
  \rput[u](24,-2){$L_2$}
  \rput[u](26,-2){$L_3$}

  \qdisk(22,-1){1.5pt}
  \qdisk(24,-1){1.5pt}
  \qdisk(26,-1){1.5pt}
  \qdisk(26,9){1.5pt}
  \qdisk(24,9){1.5pt}
  \qdisk(22,9){1.5pt}

  \rput(14.5,4){$
      \begin{array}{l}
          c_{12}'=\min(c_{12},c_{13}+c_{32})\\
          c_{13}'=c_{12}+c_{23}\\
          c_{23}'=\max(c_{23},c_{21}+c_{13}) \\[.4in]
          c_{12}=\max(c_{12}',c_{13}'+c_{32}')\\
          c_{13}=c_{12}'+c_{23}'\\
          c_{23}=\min(c_{23}',c_{21}'+c_{13}')
      \end{array}$
  }
\endpspicture
\end{center}

Enumerating rigorous paths in these graphs, we obtain
$$
\left(
\begin{array}{ccc}
M_{U_1,L_1}^{{121},3} & M_{U_1,L_2}^{{121},3} & M_{U_1,L_3}^{{121},3} \\[.1in]
M_{U_2,L_1}^{{121},3} & M_{U_2,L_2}^{{121},3} & M_{U_2,L_3}^{{121},3} \\[.1in]
M_{U_3,L_1}^{{121},3} & M_{U_3,L_2}^{{121},3} & M_{U_3,L_3}^{{121},3} 
\end{array}
\right) = 
\left(
\begin{array}{ccc}
c_{12}+c_{23}  & \min(c_{12}, c_{13}+c_{32}) & 0 \\[.1in]
+\infty        & c_{21}+c_{13}+c_{32}        & c_{21}  \\[.1in]
+\infty        & +\infty                     & c_{31}
\end{array}
\right)\,,
$$ 
\bigskip
$$
\left(
\begin{array}{ccc}
M_{U_1,L_1}^{{212},3} & M_{U_1,L_2}^{{212},3} & M_{U_1,L_3}^{{212},3} \\[.1in]
M_{U_2,L_1}^{{212},3} & M_{U_2,L_2}^{{212},3} & M_{U_2,L_3}^{{212},3} \\[.1in]
M_{U_3,L_1}^{{212},3} & M_{U_3,L_2}^{{212},3} & M_{U_3,L_3}^{{212},3} 
\end{array}
\right) = 
\left(
\begin{array}{ccc}
c_{13}' & c_{12}'                 &  0 \\[.1in]
+\infty & c_{23}'+c_{31}'+c_{12}' & \min(c_{21}',c_{23}'+c_{31}')\\[.1in]
+\infty & +\infty                 & c_{32}'+c_{21}'
\end{array}
\right)\,.
$$
It is immediate from the formulas for the transition maps $T_{121}^{212}$
and $T_{212}^{121}$ that these two matrices are equal to each other.

In the next case ($s=2$) the pseudo-lines with the lower ends $L_1$ and $L_2$
are oriented downward and the pseudo-line with the lower end $L_3$ is
oriented upward as shown on the following picture:
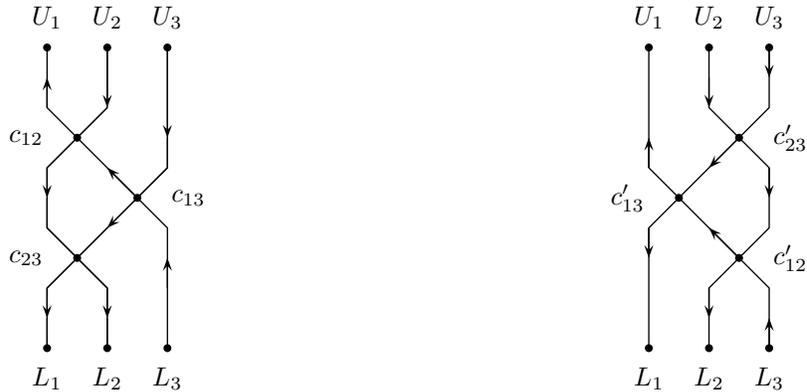
\psset{curvature=0 0 0}
\psset{xunit=.4cm,yunit=.4cm}
\begin{center}
\pspicture[.1](0,-3)(28,11)

  \psline[linewidth=.5pt, showpoints=false]{-}%
          (2,-1)(2,0)(2,1)(3,2)(4,3)(5,4)(6,5)(6,6)(6,7)(6,8)(6,9)
  \psline[linewidth=.5pt, showpoints=false]{-}%
          (4,-1)(4,0)(4,1)(3,2)(2,3)(2,4)(2,5)(3,6)(4,7)(4,8)(4,9)
  \psline[linewidth=.5pt, showpoints=false]{-}%
          (6,-1)(6,0)(6,1)(6,2)(6,3)(5,4)(4,5)(3,6)(2,7)(2,8)(2,9)

  \psline[linewidth=.5pt, showpoints=false]{-}(2,-1)(2,0)
  \psline[linewidth=.5pt, showpoints=false]{<-}(2,0)(2,1)(4,3)
  \psline[linewidth=.5pt, showpoints=false]{<-}(4,3)(6,5)(6,6)
  \psline[linewidth=.5pt, showpoints=false]{<-}(6,6)(6,9)

  \psline[linewidth=.5pt, showpoints=false]{-}(4,-1)(4,0)
  \psline[linewidth=.5pt, showpoints=false]{<-}(4,0)(4,1)(2,3)(2,4)
  \psline[linewidth=.5pt, showpoints=false]{<-}(2,4)(2,5)(4,7)
  \psline[linewidth=.5pt, showpoints=false]{<-}(4,7)(4,9)

  \psline[linewidth=.5pt, showpoints=false]{->}(6,1)(6,2)
  \psline[linewidth=.5pt, showpoints=false]{->}(5,4)(4,5)
  \psline[linewidth=.5pt, showpoints=false]{->}(2,7)(2,8)

  \qdisk(3,2){1.5pt}
  \qdisk(5,4){1.5pt}
  \qdisk(3,6){1.5pt}
  \qdisk(2,-1){1.5pt}
  \qdisk(4,-1){1.5pt}
  \qdisk(6,-1){1.5pt}
  \qdisk(6,9){1.5pt}
  \qdisk(4,9){1.5pt}
  \qdisk(2,9){1.5pt}
  \rput[u](1.3,6){$c_{12}$}
  \rput[u](6.7,4){$c_{13}$}
  \rput[u](1.3,2){$c_{23}$}
  \rput[u](2,10){$U_1$}
  \rput[u](4,10){$U_2$}
  \rput[u](6,10){$U_3$}
  \rput[u](6,-2){$L_3$}
  \rput[u](4,-2){$L_2$}
  \rput[u](2,-2){$L_1$}

  \psline[linewidth=.5pt, showpoints=false]{-}%
          (22,-1)(22,0)(22,1)(22,2)(22,3)(23,4)(24,5)(25,6)(26,7)(26,8)(26,9)
  \psline[linewidth=.5pt, showpoints=false]{-}%
          (24,-1)(24,0)(24,1)(25,2)(26,3)(26,4)(26,5)(25,6)(24,7)(24,8)(24,9)
  \psline[linewidth=.5pt, showpoints=false]{-}%
          (26,-1)(26,0)(26,1)(25,2)(24,3)(23,4)(22,5)(22,6)(22,7)(22,8)(22,9)

  \psline[linewidth=.5pt, showpoints=false]{<-}(22,2)(22,3)
  \psline[linewidth=.5pt, showpoints=false]{<-}(24,5)(25,6)
  \psline[linewidth=.5pt, showpoints=false]{<-}(26,8)(26,9)

  \psline[linewidth=.5pt, showpoints=false]{<-}(24,0)(24,1)
  \psline[linewidth=.5pt, showpoints=false]{<-}(26,4)(26,5)
  \psline[linewidth=.5pt, showpoints=false]{<-}(24,7)(24,8)

  \psline[linewidth=.5pt, showpoints=false]{->}(26,-1)(26,0)
  \psline[linewidth=.5pt, showpoints=false]{->}(25,2)(24,3)
  \psline[linewidth=.5pt, showpoints=false]{->}(22,5)(22,6)

  \qdisk(25,2){1.5pt}
  \qdisk(23,4){1.5pt}
  \qdisk(25,6){1.5pt}
  \rput[u](26.7,6){$c_{23}'$}
  \rput[u](21.3,4){$c_{13}'$}
  \rput[u](26.7,2){$c_{12}'$}
  \rput[u](22,10){$U_1$}
  \rput[u](24,10){$U_2$}
  \rput[u](26,10){$U_3$}
  \rput[u](22,-2){$L_1$}
  \rput[u](24,-2){$L_2$}
  \rput[u](26,-2){$L_3$}

  \qdisk(22,-1){1.5pt}
  \qdisk(24,-1){1.5pt}
  \qdisk(26,-1){1.5pt}
  \qdisk(26,9){1.5pt}
  \qdisk(24,9){1.5pt}
  \qdisk(22,9){1.5pt}
\endpspicture
\end{center}

In this case we have
$$
\left(
\begin{array}{ccc}
M_{L_3,U_1}^{{121},2} & M_{L_3,L_2}^{{121},2} & M_{L_3,L_1}^{{121},2} \\[.1in]
M_{U_2,U_1}^{{121},2} & M_{U_2,L_2}^{{121},2} & M_{U_2,L_1}^{{121},2} \\[.1in]
M_{U_3,U_1}^{{121},2} & M_{U_3,L_2}^{{121},2} & M_{U_3,L_1}^{{121},2} 
\end{array}
\right) = 
\left(
\begin{array}{ccc}
0       & \min(c_{12}, c_{13}+c_{32}) & c_{12}+c_{23}\\[.1in]
c_{21}  &  0        & c_{23} \\[.1in]
c_{31}  & \min(c_{32}, c_{31}+c_{12})  & c_{31}+c_{12}+c_{23}
\end{array}
\right)\,,
$$ 
\bigskip
$$
\left(
\begin{array}{ccc}
M_{L_3,U_1}^{{212},2} & M_{L_3,L_2}^{{212},2} & M_{L_3,L_1}^{{212},2} \\[.1in]
M_{U_2,U_1}^{{212},2} & M_{U_2,L_2}^{{212},2} & M_{U_2,L_1}^{{212},2} \\[.1in]
M_{U_3,U_1}^{{212},2} & M_{U_3,L_2}^{{212},2} & M_{U_3,L_1}^{{212},2} 
\end{array}
\right) = 
\left(
\begin{array}{ccc}
0       & c_{12}'                     & c_{13}' \\[.1in]
\min(c_{21}',c_{23}'+c_{31}') & 0 & \min(c_{23}',c_{21}'+c_{13}') \\[.1in]
c_{32}' +c_{21}'  & c_{32}'   & c_{32}'+c_{21}'+c_{13}' 
\end{array}
\right)\,.
$$
Again, it is clear that these two matrices are equal to each other.

The cases $s=0$ and $s=1$ are completely symmetric to the cases $s=3$ and $s=2$,
respectively. 

\medskip
We can now verify the statement of the theorem for an arbitrary~$n$.
This general statement reduces to the case of $S_3$ ($n=3$) as follows.
Clearly, it is enough to prove the statement for two reduced words $\a$ and
$\a'$ related by a $3$-move.  The corresponding reflection orderings of
inversions differ only in tree terms: $\dots < (j,k)<(i,k)<(i,j)<\dots$ and
$\dots < (i,j)<(i,k)<(j,k)<\dots$ for some $i<j<k$.  The transition map
$T_\a^{\a'}$ is the map $T_{ijk}$ that transforms $c_{ij}$, $c_{ik}$, and
$c_{jk}$ into $c_{ij}'$, $c_{ik}'$, and $c_{jk}'$ according
formulas~(\ref{eq:transform_general}) and does not change other variables.

The intersection points of the pseudo-lines labelled $i$, $j$, and~$k$ form a
subdiagram $S$ in the wiring diagram of $\a$ (respectively, 
a subdiagram $S'$ in the wiring diagram of $\a'$) isomorphic to a
wiring diagram for $S_3$.  Let us add six auxiliary vertices $u_1$, $u_2$,
$u_3$, and $l_1$, $l_2$, $l_3$ to the graph  $G(\a,s)$ (respectively, in
$G(\a',s)$) that mark the upper and the lower ends of the pseudo-lines $i$, $j$, and
$k$ in this subdiagram.

If a path $P$ in the graph $G(\a,s)$ does not pass through any of the 
vertices $v_{ij}$, $v_{ik}$, and $v_{jk}$ then the expression~(\ref{eq:c_P}) 
for $c_P$ does not change under the transformation map $T_\a^{\a'}$.
Otherwise, the path $P$ arrives to the subdiagram $S$ via one of the
six auxiliary points $u_1$,\dots,$l_3$ and leaves the subdiagram via
another of these six points.  

Let us fix two vertices $b$ and $e$ of the six auxiliary vertices 
and two rigorous paths $P_1$ (from $B$ to $b$) and $P_2$ (from $e$ to $E$).
And let $\bar M_{b,e,P_1,P_2}^{\a,s}(C)$ 
(respectively, $\bar M_{b,e,P_1,P_2}^{\a',s}(C')$)
be the minimum of the expressions $c_P$ over rigorous paths $P$ in $G(\a,s)$ 
(respectively, in $G(\a',s)$) which are obtained by concatenation of
the path $P_1$, a rigorous path in $S$ (respectively, in $S'$) 
from $b$ to $e$, and the path~$P_2$.
Then, by our definitions, 
$$
M_{B,E}^{\a,s}=\min_{b,e,P_1,P_2} \bar M_{b,e,P_1,P_2}^{\a,s}
\quad
\textrm{and}
\quad 
M_{B,E,s}^{\a',s}=\min_{b,e,P_1,P_2} \bar M_{b,e,P_1,P_2}^{\a',s}.
$$
It follows from the case of $S_3$ considered above that 
$\bar M^\a_{b,e,P_1,P_2}(C)= \bar M^{\a'}_{b,e,P_1,P_2}(C')$.
Therefore, $M_{B,E}^{\a,s}(C)=M_{B,E}^{\a',s}(C')$. 
This proves the theorem.
\endproof

\medskip
\proofof{Theorem~\ref{th:low}}
Suppose that $low(\b)=(i,j)$, the lower end of $i$-th pseudo-line is $L_{s+1}$,
and the lower end of $j$-th pseudo-line is $L_s$.  (In the case of $\wnot$ 
we have $i=n-s$ and $j=n-s+1$.) 
Then there is only one rigorous path in $G(\b,s)$ 
from $L_{s+1}$ to $L_s$, namely, $L_{s+1}\to v_{ij}\to L_s$.
In this case $M_{L_{s+1},L_s}^{\b,s}(C')= c_{ij}'$.
Thus $M_{L_{s+1},L_s}^{\a,s}(C)=M_{L_{s+1},L_s}^{\b,s}(C')=
c_{ij}'=c_{low(\b)}'$, where $C'=T_\a^\b(C)$.

For any $s=1,\dots,n-1$, there is a reduced decomposition $\b$ of the longest
permutation $\wnot$ 
such that $low(\b)=(n-s,n-s+1)$.
Thus the inequality $M_{L_{s+1},L_s}^{\a,s}(C)\geq 0$ is equivalent
to saying that for any reduced word $\b\in\R(\wnot)$ such that $low(\b)=(n-s,n-s+1)$
the lowest entry $c'_{low(\b)}$ of $C'=T_\a^\b(C)$ is nonnegative.
The statement follows.
\endproof

\medskip
\proofof{Theorem~\ref{th:princ-w0}}  
We deduce this theorem from Theorem~\ref{th:low}.
Let us fix $\a\in\R(\wnot)$.
It is enough to show that if $C=(c_{ij})$ has a negative entry then
there is a reduced word $\b\in\R(\wnot)$ such that the lowest entry $c'_{low(\b)}$
of $C'=T_\a^\b(C)$ is negative.

Suppose not. 
Let us pick a reduced word $\b$ such that $C'=T_\a^{\b}(C)$ 
has a negative entry $c_{pq}'<0$ located on the lowest possible level.
The pair $(p,q)\ne low(\b)$ does not correspond to the lowest crossing
in the wiring diagram of $\b$. 
Thus (possibly, after several 3-moves that don't 
affect $c_{pq}'$) we can make a 3-move transforming three entries
$(c_{ij}',c_{ik}',c_{jk}')\to (c_{ij}'',c_{ik}'',c_{jk}'')$
by the rule~(\ref{eq:transform_general}) such that 
$(p,q)\in\{(i,j),(i,k),(j,k)\}$
but $(p,q)$ is not the lowest pair $(j,k)$ among these three.
By our assumption, $c_{jk}'$ is nonnegative.
Then $c_{ij}''=\min(c_{ij}', c_{ik}'-c_{jk})$ is negative
and $c_{ij}''$ is located on a lower level in the resulting 
wiring diagram than the level of $c_{pq}'$ in $\b$.
Contradiction.
\endproof

\section{Associativity}                         
\label{sec:associativity}
\neweq

In this section we prove Theorem~\ref{th:associative}, which claims
that the $*$-product defined by~(\ref{eq:mmn}) 
is an associative operation.  

\proofof{Theorem~\ref{th:associative}}
We need to verify that
\begin{equation}
(e_{x_1\dots x_m} * e_{y_1\dots y_n})* e_{z_1\dots z_k} =
e_{x_1\dots x_m} * (e_{y_1\dots y_n}* e_{z_1\dots z_k}),
\label{eq:exyz}
\end{equation}
for any positive $m,n,k$ and 
$x_1\leq \cdots\leq  x_m$, $y_1\leq \cdots\leq y_n$,
$z_1\leq \cdots \leq z_k$.

Let $\Id_k$ be the identity permutation in $S_k$.
The permutation $w(m,n)\times \Id_k\in S_{m+n}\times S_k$ 
is canonically embedded into $S_{m+n+k}$.
Likewise, the permutation $\Id_m\times w(n,k)\in S_m\times S_{n+k}$ 
is canonically embedded into $S_{m+n+k}$.
Then $w(m+n,k)\cdot \left(w(m,n)\times \Id_k\right) =
w(m,n+k)\cdot \left(\Id_m\times w(n,k)\right)$.
We will denote this permutation by $w(m,n,k)$. 

Remind that the permutations $w(m+n,k)$ and $w(m,n)\times \Id_k$ have
unique (up to $2$-moves) reduces decompositions.  Let $\a^1$ be
a reduced word for $w(m,n,k)$ obtained by concatenation of 
reduced words for $w(m+n,k)$ and $w(m,n)\times \Id_k$.
Analogously, let $\a^2$ be a reduced word for $w(m,n,k)$ obtained by 
concatenation of reduced words for $w(m,n+k)$ and $\Id_m\times w(n,k)$.

The inversion set $I(w(m,n,k))$ of the permutation $w(m,n,k)$ is
the union of the following three sets of pairs 
$[1,m]\times [m+1,m+n]$, 
$[1,m]\times [m+n+1,m+n+k]$, and
$[m+1,m+n]\times [m+n+1,m+n+k]$,
where $[a,b]=\{a,a+1,\dots,b\}$.

By the definition of $*$-product, the left hand side of the 
expression~(\ref{eq:exyz}) is equal to the sum
$\sum R_{\a^1}(C) \cdot (e_{x_1\dots x_m} \otimes e_{y_1\dots y_n}
\otimes e_{z_1\dots z_k})$
over all collections $C=(c_{ij}) \in\Z^{I(w(m,n,k))}$ with nonnegative
integer entries such that 
$$
\begin{array}{l}
c_{ij}\geq c_{pq}\qquad\textrm{whenever}
\quad 1\leq p\leq i\leq m,\quad m+1 \leq j\leq q\leq m+n,  \\[.05in]
c_{ij}\geq c_{pq}\qquad\textrm{whenever}
\quad 1\leq p\leq i\leq m+n,\quad m+n+1\leq j\leq q\leq m+n+k, 
\end{array}
$$
cf.~(\ref{eq:conditions}).
These are exactly the inequalities defining the principal cones
$\C_{\a^1}$, cf.~Example~\ref{ex:rmn}.  
Thus the left hand side of~(\ref{eq:exyz}) can be written as 
$$
\sum_{C\in\C_{\a^1}}  R_{\a^1}(C)
\cdot (e_{x_1\dots x_m}\otimes e_{y_1\dots y_n}\otimes e_{z_1\dots z_k}).
$$
Analogously, the right hand side of~(\ref{eq:exyz}) can be written as 
$$
\sum_{C\in\C_{\a^2}} R_{\a^2}(C)
\cdot (e_{x_1\dots x_m} \otimes e_{y_1\dots y_n} \otimes e_{z_1\dots z_k}).
$$
The equality of these two expressions follows from~(\ref{eq:Rab}) and
Theorem~\ref{th:cones-bijectively}.

\medskip
This proves Theorem~\ref{th:associative}
and thus completes the proof of our main statement concerning 
the $*$-product---Theorem~\ref{th:main}.
\endproof

\section{Web Functions, Berenstein-Zelevinsky Triangles,\\ and Hidden Duality}
\label{sec:web}

In this section we give a geometric interpretation of the scattering 
matrix~(\ref{eq:scattering}) in terms of certain {\it web functions}
as well as a ``physical'' motivation for it.
Then we establish a relationship between integral web diagrams and fillings of 
Berenstein-Zelevinsky triangles.  We also discuss the ``hidden duality'' of 
the Littlewood-Richardson coefficients under conjugation of partitions:
$c_{\lambda \mu}^\nu = c_{\lambda' \mu'}^{\nu'}$.

\subsection{Web functions}

It will be convenient to use the baricentric coordinates in~$\RR^2$.
Namely, we represent a point in~$\RR^2$ by a triple $(\alpha,\beta,\gamma)$
such that $\alpha+\beta+\gamma=0$.  
We say that a line in~$\RR^2$ is of the {\it first}
(respectively, {\it second,} or {\it third\/}) {\it type} if its 
first (respectively, second, or third) baricentric coordinate is fixed.
We will denote by $(a,*,*)$ the first type line given by 
$\{(\alpha,\beta,\gamma)\mid \alpha+\beta+\gamma=0, \ \alpha = a\}$.  
Analogously, we will denote by $(*,b,*)$ and $(*,*,c)$ the lines of the second 
and third types given by
$\{(\alpha,\beta,\gamma)\mid \alpha+\beta+\gamma=0, \ \beta = b\}$ and  
$\{(\alpha,\beta,\gamma)\mid \alpha+\beta+\gamma=0, \ \gamma = c\}$,
respectively.  
Each of the following two pictures represents a union of three rays 
of first, second, and third type originating at the same point.
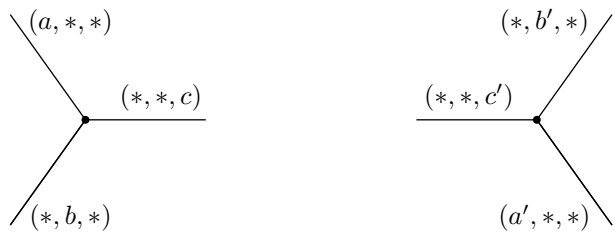
\psset{curvature=0 0 0}
\psset{xunit=.2cm,yunit=.2cm}
\begin{center}
\pspicture[.1](0,-1)(40,15)
  \pscurve[linewidth=.5pt, showpoints=false]{-}%
          (0,0)(5,7)(13,7)
  \pscurve[linewidth=.5pt, showpoints=false]{-}%
          (0,0)(5,7)(0,14)
  \qdisk(5,7){1.5pt}
  \rput[u](4,13.5){$(a,*,*)$}
  \rput[u](4,0.5){$(*,b,*)$}
  \rput[u](10,8.5){$(*,*,c)$}

  \pscurve[linewidth=.5pt, showpoints=false]{-}%
          (40,0)(35,7)(27,7)
  \pscurve[linewidth=.5pt, showpoints=false]{-}%
          (40,0)(35,7)(40,14)
  \qdisk(35,7){1.5pt}
  \rput[u](35.5, 0.5){$(a',*,*)$}
  \rput[u](35.5, 13.5){$(*,b',*)$}
  \rput[u](30.5, 8.5){$(*,*,c')$}
\endpspicture
\end{center}
Notice that in both cases we have $a+b+c=0$ and $a'+b'+c'=0$.
We say that these two types of sets are {\it left} and {\it right forks.}
The central point of a fork is called its {\it node.}  The node
of the left (respectively, right) fork shown above is the point~$(a,b,c)$
(respectively, $(a',b',c')$),  in the baricentric coordinates.
We say that a function $f:\RR^2\to \RR$ is a {\it fork function} (left or
right) if there is a fork such that $f$ is equal to~1 on three rays of the fork,
to~$3/2$ on its node, and~$0$ everywhere else.

\begin{definition} {\rm
A {\it web function} is a function 
$f:\RR^2\to \RR_+$
such that for every point in~$\RR^2$ there exists an open neighborhood~$U$ of the point,
for which the restriction~$f|_U$ is either zero, or the characteristic
function of a line of one of three types, or a fork function (left or right), 
or a finite sum of several such functions. 
We say that a web function is {\it integral} if all its lines are
of the form $(a,*,*)$, $(*,b,*)$, or~$(*,*,c)$ with integer~$a$,
$b$, and~$c$.
}
\end{definition}

Geometrically, we represent a web function by a picture (called {\it web diagram\/}) 
composed of rays and line segments of one of three types, 
and left or right forks (possibly doubled, tripled, etc.).
See below for examples of web diagrams.

We say that a web function is {\it generic} if it only takes values~$0$, $1$
and~$3/2$. In other words, the diagram of a generic web function is composed 
of {\it noncrossing} rays, line segments, and forks.  
An arbitrary web diagram can be obtained by degeneration 
of a generic web diagram, i.e., by merging several lines, line intervals, and nodes 
together.  For example the following diagram on the left-hand side presents a generic web
function.  The diagram of a non-generic web function on the 
right-hand side is obtained by merging three nodes together.  The
double line shows the locus where the web function is equal to~$2$.
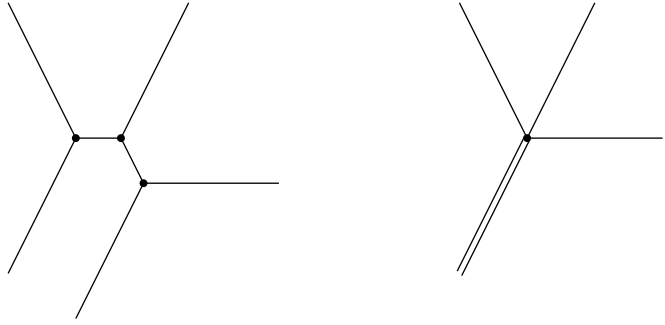
\psset{xunit=.3cm,yunit=.3cm}
\begin{center}
\pspicture[.1](-1,-3)(31,13)
  \psline[linewidth=.5pt]{-}(-1,0)(2,6)(-1,12)
  \psline[linewidth=.5pt]{-}(2,6)(4,6)
  \psline[linewidth=.5pt]{-}(2,-2)(5,4)(4,6)(7,12)
  \psline[linewidth=.5pt]{-}(5,4)(11,4)
  \qdisk(2,6){1.5pt}
  \qdisk(5,4){1.5pt}
  \qdisk(4,6){1.5pt}

  \psline[linewidth=.5pt]{-}(19.1, -0.1)(22.1, 5.9)
  \psline[linewidth=.5pt]{-}(18.9, 0.1)(21.9, 6.1)
  \psline[linewidth=.5pt]{-}(22,6)(19,12)
  \psline[linewidth=.5pt]{-}(22,6)(28,6)
  \psline[linewidth=.5pt]{-}(22,6)(25,12)
  \qdisk(22,6){1.5pt}
\endpspicture
\end{center}
Each web diagram consists of several nodes, line intervals,
semi-infinite rays, and/or infinite lines.  We are only interested
in web functions whose diagrams have finitely many nodes.  We will refer to 
semi-infinite rays in a web diagram as {\it boundary rays.}
It is possible that the boundary rays are doubled (as in the example
above), tripled, etc.  The possible directions for boundary rays are
{\it North-West} and {\it South-East} (for type~1 rays), 
{\it North-East} and {\it South-West} (for type~2 rays), and
{\it West} and {\it East} (for type~3 rays), as shown on the picture 
below:
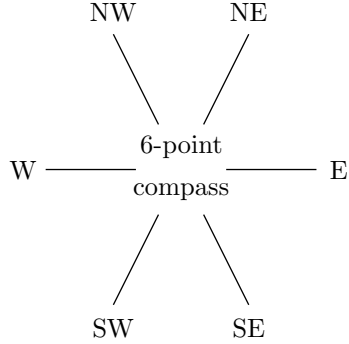
\psset{xunit=.3cm,yunit=.3cm}
\begin{center}
\pspicture[.1](-3,-2)(9,14)
  \psline[linewidth=.5pt]{-}(0,0)(2,4)
  \psline[linewidth=.5pt]{-}(6,0)(4,4)
  \psline[linewidth=.5pt]{-}(-3,6)(1,6)
  \psline[linewidth=.5pt]{-}(5,6)(9,6)
  \psline[linewidth=.5pt]{-}(0,12)(2,8)
  \psline[linewidth=.5pt]{-}(6,12)(4,8)
  \rput[u](0,-1){SW}
  \rput[u](6,-1){SE}
  \rput[u](0,13){NW}
  \rput[u](6,13){NE}
  \rput[u](-4,6){W}
  \rput[u](10,6){E}
  \rput[u](3,7){$6$-point}
  \rput[u](3,5){compass}
\endpspicture
\end{center}

Recall that we defined the scattering matrix~$R(c)$ by
$$
R(c)\,:\, e_x\otimes e_y \longmapsto 
\left\{
\begin{array}{cl}
e_{y+c}\otimes e_{x-c} & \textrm{if } c \geq x-y, \\[.1in]
0                      & \textrm{otherwise,}
\end{array}
\right.
$$
(see Definition~\ref{def:scattering}).  Here 
(unlike Section~\ref{sec:scattering}) we allow 
$x$, $y$, and~$c$ to be any real numbers.

The first type line $(-x,*, *)$ can be thought of as the trajectory
of a certain {\it left particle} of energy~$x$.  
We will denote this particle by $\rmleft(x)$.
Analogously, the second type line $(*,y,*)$ represents the trajectory of 
a {\it right particle} of energy~$y$, denoted by $\rmright(y)$.
In both cases the trajectories go downward (from left to right for left
particles and from right to left for right particles).  Then the scattering 
matrix~$R(c)$ represents an interaction of a left particle 
of energy~$x$ with a right particle
of energy~$y$.  The following web diagram
visualizes the scattering matrix~$R(c)$.
\psset{xunit=.3cm,yunit=.3cm}
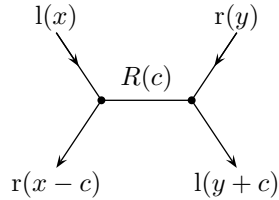
\begin{center}
\pspicture[.1](0,0)(10,8)
  \psline[linewidth=.5pt]{<-}(2,1)(4,4)(2,7)
  \psline[linewidth=.5pt]{<-}(3,5.5)(2,7)
  \psline[linewidth=.5pt]{<-}(10,1)(8,4)(10,7)
  \psline[linewidth=.5pt]{<-}(9,5.5)(10,7)
  \psline[linewidth=.5pt,showpoints=false]{-}(4,4)(8,4)
  \qdisk(4,4){1.5pt}
  \qdisk(8,4){1.5pt}
  \rput[u](2,.3){$\rmright(x-c)$}
  \rput[u](10,.3){$\rmleft(y+c)$}
  \rput[u](2,7.7){$\rmleft(x)$}
  \rput[u](10,7.7){$\rmright(y)$}
  \rput[u](6,5){$R(c)$}
\endpspicture
\end{center}
The horizontal segment in this diagram lies on the third type line $(*,*,c)$.
Thus the interaction~$R(c)$ happens on the level~$c$.  The condition
$c\geq x-y$ means that the interaction happens before the trajectories
of the particles $\rmleft(x)$ and $\rmright(y)$ cross each other.

Recall that in Section~\ref{sec:scattering}
we defined the operator~$R_{(m,n)}\((c_{ij})\)$ as a composition of 
of the scattering matrices~$R_{ij}(c_{ij})$, $1\leq i\leq m$, 
$m+1\leq j\leq m+n$.
The operator~$R_{(m,n)}\((c_{ij})\)$ applied to the 
vector~$e_{x_1}\otimes \cdots \otimes e_{x_m}\otimes
e_{y_1}\otimes\cdots\otimes e_{y_n}$ and producing the vector
$e_{z_1}\otimes \cdots \otimes e_{z_{n+m}}$ can be represented by a
web digram, which is a combinations of several pieces similar
to the one shown above.  
In our pseudophysical lexicon,
this diagram represents an interaction of $m$ left particles
with $n$ right particles.
An example of such a web diagram 
for $m=4$ and $n=3$ is given below:
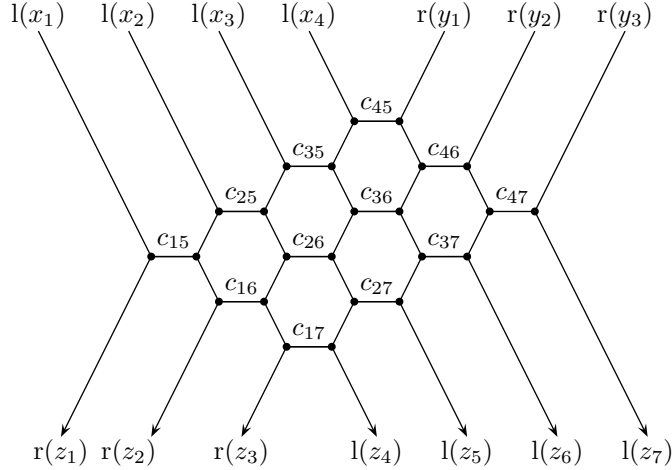
\psset{xunit=.3cm,yunit=.3cm}
\begin{center}
\pspicture[.1](-5,-4)(20,18)
  \psline[linewidth=.5pt]{<-}(-5,-2)(-1,6)(1,6) (2,8)(-2,16)
  \psline[linewidth=.5pt]{<-}(-1,-2)(2,4)(4,4) (5,6)(4,8)(5,10)(2,16)
  \psline[linewidth=.5pt]{<-}(3,-2)(5,2)(7,2) (8,4) (7,6) (8,8)(7,10)(8,12)(6,16)

  \psline[linewidth=.5pt]{<-}(13,-2)(10,4)(11,6)(10,8)(11,10)(10,12)(12,16)
  \psline[linewidth=.5pt]{<-}(17,-2)(13,6)(14,8)(13,10)(16,16)
  \psline[linewidth=.5pt]{<-}(21,-2)(16,8)(20,16)

  \psline[linewidth=.5pt]{<-}(9,-2)(7,2)(5,2)(4,4)(2,4)(1,6)(-1,6)(-6,16)

  \psline[linewidth=.5pt]{-}(2,8)(4,8)
  \psline[linewidth=.5pt]{-}(2,8)(4,8)
  \psline[linewidth=.5pt]{-}(5,6)(7,6)
  \psline[linewidth=.5pt]{-}(5,10)(7,10)
  \psline[linewidth=.5pt]{-}(8,4)(10,4)
  \psline[linewidth=.5pt]{-}(8,8)(10,8)
  \psline[linewidth=.5pt]{-}(8,12)(10,12)
  \psline[linewidth=.5pt]{-}(11,6)(13,6)
  \psline[linewidth=.5pt]{-}(11,10)(13,10)

  \psline[linewidth=.5pt]{-}(5,6)(7,6)
  \psline[linewidth=.5pt]{-}(5,10)(7,10)
  \psline[linewidth=.5pt]{-}(8,4)(10,4)
  \psline[linewidth=.5pt]{-}(8,8)(10,8)
  \psline[linewidth=.5pt]{-}(8,12)(10,12)
  \psline[linewidth=.5pt]{-}(11,6)(13,6)
  \psline[linewidth=.5pt]{-}(11,10)(13,10)
  \psline[linewidth=.5pt]{-}(14,8)(16,8)

  \qdisk(7,2){1.5pt}
  \qdisk(5,2){1.5pt}
  \qdisk(4,4){1.5pt}
  \qdisk(2,4){1.5pt}
  \qdisk(1,6){1.5pt}
  \qdisk(-1,6){1.5pt}

  \qdisk(2,8){1.5pt}

  \qdisk(5,6){1.5pt}
  \qdisk(4,8){1.5pt}
  \qdisk(5,10){1.5pt}

  \qdisk(8,4){1.5pt}
  \qdisk(7,6){1.5pt}
  \qdisk(8,8){1.5pt}
  \qdisk(7,10){1.5pt}
  \qdisk(8,12){1.5pt}

  \qdisk(10,4){1.5pt}
  \qdisk(11,6){1.5pt}
  \qdisk(10,8){1.5pt}
  \qdisk(11,10){1.5pt}
  \qdisk(10,12){1.5pt}

  \qdisk(13,6){1.5pt}
  \qdisk(14,8){1.5pt}
  \qdisk(13,10){1.5pt}

  \qdisk(16,8){1.5pt}

  \rput[u](-6,16.7){$\rmleft(x_1)$}
  \rput[u](-2,16.7){$\rmleft(x_2)$}
  \rput[u](2,16.7){$\rmleft(x_3)$}
  \rput[u](6,16.7){$\rmleft(x_4)$}
  \rput[u](12,16.7){$\rmright(y_1)$}
  \rput[u](16,16.7){$\rmright(y_2)$}
  \rput[u](20,16.7){$\rmright(y_3)$}

  \rput[u](-5,-2.7){$\rmright(z_1)$}
  \rput[u](-2,-2.7){$\rmright(z_2)$}
  \rput[u](3,-2.7){$\rmright(z_3)$}
  \rput[u](9,-2.7){$\rmleft(z_4)$}
  \rput[u](13,-2.7){$\rmleft(z_5)$}
  \rput[u](17,-2.7){$\rmleft(z_6)$}
  \rput[u](21,-2.7){$\rmleft(z_7)$}

  \rput[u](3,8.7){$c_{25}$}
  \rput[u](6,6.7){$c_{26}$}
  \rput[u](6,10.7){$c_{35}$}
  \rput[u](9,4.7){$c_{27}$}
  \rput[u](9,8.7){$c_{36}$}
  \rput[u](9,12.7){$c_{45}$}
  \rput[u](12,6.7){$c_{37}$}
  \rput[u](12,10.7){$c_{46}$}
  \rput[u](15,8.7){$c_{47}$}

  \rput[u](6,2.7){$c_{17}$}
  \rput[u](3,4.7){$c_{16}$}
  \rput[u](0,6.7){$c_{15}$}
\endpspicture
\end{center}
In general, such a web diagram need not be 
as regular as the one shown above.  
The edge lengths can be arbitrarily deformed.  

This web diagram has the following boundary rays: North-West rays
corresponding to incoming particles $\rmleft(x_1),\dots,\rmleft(x_m)$;
North-East rays, corresponding to incoming particles
$\rmright(y_1),\dots,\rmright(y_n)$; South-West rays, corresponding to outgoing
particles $\rmright(z_1),\dots,\rmright(z_n)$; South-East rays, corresponding
to outgoing particles $\rmleft(z_{n+1}),\dots,\rmleft(z_{n+m})$; and  no East
or West boundary rays.  The $i$-th left patricle interacts with the $j$-th
right particle on the level $c_{i\,j+m}$.  In the web diagram, this interaction
is represented by an interval which lies on the line $(*,*,c_{i\,j+m})$.
Such a web digram is integral if and only if all
$x_i$, $y_j$, $z_k$, and $c_{ij}$ are integers.

Using the language of web diagrams, we derive the following statement
from Corollary~\ref{cor:rule}.

\begin{corollary}
\label{cor:web}
Let $\lambda = \omega_{x_1}+\cdots +\omega_{x_m}$,
$\mu = \omega_{y_1}+\cdots +\omega_{y_n}$, and
$\nu = \omega_{z_1}+\cdots +\omega_{z_{m+n}}$ be three dominant weights in
$GL(N)$, where 
$1\leq x_1\leq \cdots \leq x_m\leq N$, 
$1\leq y_1\leq \cdots \leq y_n\leq N$, and
$0\leq z_1\leq \cdots \leq z_{m+n}\leq N$. 
The Littlewood-Richardson coefficient $c_{\lambda\mu}^\nu$ 
is equal to the number of integral web diagrams which have the following 
fixed boundary rays:
\begin{itemize}
\item
the North-West rays $(-x_1,*,*)$, \dots, $(-x_m,*,*)$;
\item
the North-East rays $(*,y_1,*)$, \dots, $(*,y_n,*)$;
\item 
the South-West rays $(*,z_1,*)$, \dots, $(*,z_n,*)$;
\item
the South-East rays $(-z_{n+1},*,*)$, \dots, $(-z_{n+m},*,*)$;
\item
no East or West boundary rays.
\end{itemize}
\end{corollary}

Independently of our work a similar to a web diagram 
notion of a {\it honeycomb tinkertoy} recently 
appeared in~\cite{Knutson} in relation to Klyachko's saturation hypothesis.  
(The origin of the term ``honeycomb'' should be clear from the previous
picture.)  This tinkertoy is given along with a statement 
reminiscent of Corollary~\ref{cor:web}.
In our notation this statement can be reformulated as follows:

\begin{theorem} 
{\rm \cite[Theorem~1]{Knutson}} \ 
\label{th:web-dual}
Let $\lambda=(\lambda_1,\dots,\lambda_N)$,
$\mu=(\mu_1,\dots,\mu_N)$, and 
$\nu=(\nu_1,\dots,\nu_N)$ be three dominant weights of~$GL(N)$.
The Littlewood-Richardson coefficient~$c_{\lambda\mu}^\nu$ is 
equal to the number of integral web diagrams with the following 
fixed boundary rays:
\begin{itemize}
\item
the North-West rays $(\lambda_1,*,*)$,\dots,$(\lambda_N,*,*)$;
\item
the South-West rays $(*,\mu_1,*)$,\dots,$(*,\mu_N,*)$;
\item
the East rays $(*,*,-\nu_N)$,\dots,$(*,*,-\nu_1)$.
\end{itemize}
\end{theorem}

In a sense, these two statements are dual to each other.
The proof of Theorem~\ref{th:web-dual} is based on a simple one-to-one 
correspondence (see~\cite[Appendix]{Knutson}) between integral honeycomb tinkertoys 
(in our notation, web diagrams satisfying the conditions of Theorem~\ref{th:web-dual}) 
and Berenstein-Zelevinsky patterns~\cite{BZ1}.  
This correspondence just assigns to such a web diagram the triangular array filled 
by lengths of edges of the diagram.

The Berenstein-Zelevinsky interpretation of the Littlewood-Richardson coefficients, 
among its many other virtues, makes it clear that these coefficients are symmetric with 
respect to the action of~$S_3$ by permuting the three weights.  Nevertheless this
construction obscures the invariance of the Littlewood-Richardson coefficients 
under the conjugation of partitions: $c_{\lambda \mu}^\nu = 
c_{\lambda' \mu'}^{\nu'}$.  This ``hidden duality'' 
can be observed from another even simplier bijection between web diagrams and 
Berenstein-Zelevinsky patterns, which is ``dual'' to the one given 
in~\cite[Appendix]{Knutson}.  
To formulate the correspondence we have to rigorously
define these patterns.

\subsection{BZ-functions and BZ-triangles}

We say that {\it BZ-lattice} $\L_{BZ}$ is the set $({1\over 2}\Z \times {1\over 2}\Z)
\setminus (\Z\times\Z)$.  Using the baricentric coordinates we can describe $\L_{BZ}$
as the set of points $(\alpha, \beta, \gamma)$, $\alpha+\beta+\gamma=0$, such that
$2\alpha$, $2\beta$, and $2\gamma$ are integer but at least one $\alpha$,
$\beta$, or $\gamma$ is not integer.

Every integer point $(a,b,c)$, $a+b+c=0$, has six neighbours
in $\L_{BZ}$ that form the vertices of the following hexagon:
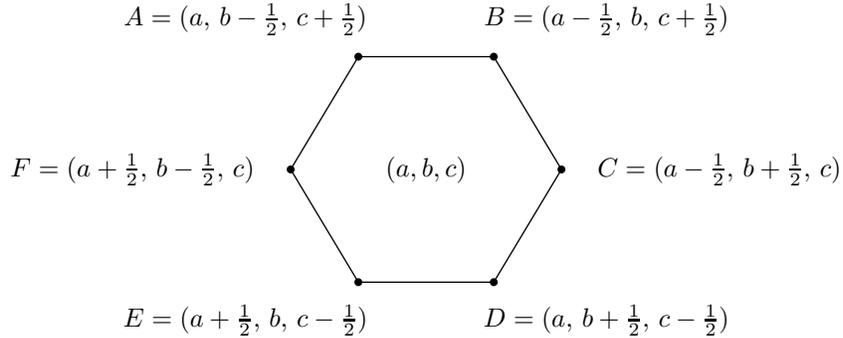
\psset{xunit=.3cm,yunit=.25cm}
\begin{center}
\pspicture[.1](-3,-4)(9,16)
  \psline[linewidth=.5pt,showpoints=false]{-}(0,0)(6,0)(9,6)(6,12)(0,12)(-3,6)(0,0)
  \qdisk(0,0){1.5pt}
  \qdisk(6,0){1.5pt}
  \qdisk(-3,6){1.5pt}
  \qdisk(9,6){1.5pt}
  \qdisk(0,12){1.5pt}
  \qdisk(6,12){1.5pt}

 \rput[u](3,6){$(a,b,c)$}
 \rput[u](-5,-2){$E=(a+{1\over 2},\,b,\,c-{1\over 2})$}
 \rput[u](11,-2){$D=(a,\,b+{1\over 2},\,c-{1\over 2})$}
 \rput[u](-5,14){$A=(a,\,b-{1\over 2},\,c+{1\over 2})$}
 \rput[u](11,14){$B=(a-{1\over 2},\,b,\,c+{1\over 2})$}
 \rput[u](-10,6){$F=(a+{1\over 2},\,b-{1\over 2},\,c)$}
 \rput[u](16,6){$C=(a-{1\over 2},\,b+{1\over 2},\,c)$}
\endpspicture
\end{center}

\begin{definition} {\rm
A function $f:\L_{BZ}\to\{0,1,2,\dots\}$ is called a 
{\it BZ-function} if for any hexagon as above
it satisfies the following {\it hexagon condition:}
$$f(A)+f(B)=f(D)+f(E),\quad f(B)+f(C)=f(E)+f(F),\textrm{ and }f(C)+f(D)=f(F)+f(A).$$
}
\end{definition}

\begin{proposition}
Integral web functions are in one-to-one correspondence with
BZ-functions.  This correspondence $\kappa$ is given by restricting a web function
$f:\RR^2\to\RR_+$ to the BZ-lattice $\L_{BZ}$:
$$
\begin{array}{rcccc}
\kappa &:& \{\textrm{integral web functions}\}
&\longrightarrow &\{\textrm{BZ-functions}\}
\\[.1in]
\kappa&:& f &\longmapsto& f\,|_{\L_{BZ}}\,.
\end{array}
$$
\label{prop:web=BZ}
\end{proposition}

\proof
Restrictions of first, second, and third type lines and of left and right forks 
to a hexagon are the following:
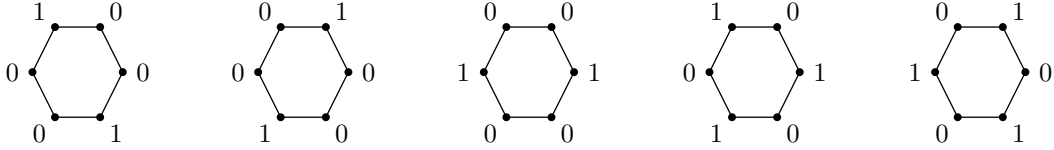
\psset{xunit=.3cm,yunit=.3cm}
\begin{center}
\pspicture[.1](-2,-1)(44,5)
  \psline[linewidth=.5pt,showpoints=false]{-}(0,0)(2,0)(3,2)(2,4)(0,4)(-1,2)(0,0)
  \qdisk(0,0){1.5pt}
  \qdisk(2,0){1.5pt}
  \qdisk(3,2){1.5pt}
  \qdisk(2,4){1.5pt}
  \qdisk(0,4){1.5pt}
  \qdisk(-1,2){1.5pt}

 \rput[u](-.7,-.7){$0$}
 \rput[u](2.7,-.7){$1$}
 \rput[u](3.9,2){$0$}
 \rput[u](2.7,4.7){$0$}
 \rput[u](-.7,4.7){$1$}
 \rput[u](-1.9,2){$0$}

  \psline[linewidth=.5pt,showpoints=false]{-}(10,0)(12,0)(13,2)(12,4)(10,4)(9,2)(10,0)
  \qdisk(10,0){1.5pt}
  \qdisk(12,0){1.5pt}
  \qdisk(13,2){1.5pt}
  \qdisk(12,4){1.5pt}
  \qdisk(10,4){1.5pt}
  \qdisk(9,2){1.5pt}

 \rput[u](9.3,-.7){$1$}
 \rput[u](12.7,-.7){$0$}
 \rput[u](13.9,2){$0$}
 \rput[u](12.7,4.7){$1$}
 \rput[u](9.3,4.7){$0$}
 \rput[u](8.1,2){$0$}

  \psline[linewidth=.5pt,showpoints=false]{-}(20,0)(22,0)(23,2)(22,4)(20,4)(19,2)(20,0)
  \qdisk(20,0){1.5pt}
  \qdisk(22,0){1.5pt}
  \qdisk(23,2){1.5pt}
  \qdisk(22,4){1.5pt}
  \qdisk(20,4){1.5pt}
  \qdisk(19,2){1.5pt}

 \rput[u](19.3,-.7){$0$}
 \rput[u](22.7,-.7){$0$}
 \rput[u](23.9,2){$1$}
 \rput[u](22.7,4.7){$0$}
 \rput[u](19.3,4.7){$0$}
 \rput[u](18.1,2){$1$}

  \psline[linewidth=.5pt,showpoints=false]{-}(30,0)(32,0)(33,2)(32,4)(30,4)(29,2)(30,0)
  \qdisk(30,0){1.5pt}
  \qdisk(32,0){1.5pt}
  \qdisk(33,2){1.5pt}
  \qdisk(32,4){1.5pt}
  \qdisk(30,4){1.5pt}
  \qdisk(29,2){1.5pt}

 \rput[u](29.3,-.7){$1$}
 \rput[u](32.7,-.7){$0$}
 \rput[u](33.9,2){$1$}
 \rput[u](32.7,4.7){$0$}
 \rput[u](29.3,4.7){$1$}
 \rput[u](28.1,2){$0$}

  \psline[linewidth=.5pt,showpoints=false]{-}(40,0)(42,0)(43,2)(42,4)(40,4)(39,2)(40,0)
  \qdisk(40,0){1.5pt}
  \qdisk(42,0){1.5pt}
  \qdisk(43,2){1.5pt}
  \qdisk(42,4){1.5pt}
  \qdisk(40,4){1.5pt}
  \qdisk(39,2){1.5pt}

 \rput[u](39.3,-.7){$0$}
 \rput[u](42.7,-.7){$1$}
 \rput[u](43.9,2){$0$}
 \rput[u](42.7,4.7){$1$}
 \rput[u](39.3,4.7){$0$}
 \rput[u](38.1,2){$1$}

\endpspicture
\end{center}
It is clear that all these five functions satisfy the hexagon condition.
It is also not hard to verify that any nonnegative integer function on a hexagon
that satisfies the hexagon condition is a linear combination of these five functions
with nonnegative integer coefficients. 
Thus restrictions of integral web functions to the BZ-lattice are BZ-functions
and every BZ-function can be obtained in such a way.  On the other hand,
an integral web function is determined by its values on $\L_{BZ}$.
For example, the values in the center of the hexagon are equal to $1$ 
for the first three functions above and to $3/2$ for the remaining two functions.
\endproof

Let us fix an integer~$N\geq 1$.  The BZ-triangle $T_N$ is the triangular subset 
in~$\L_{BZ}$ given by the inequalities $\alpha>-N$, $\beta>0$,
and $\gamma=-\alpha-\beta>0$. 
A {\it Berenstein-Zelevinsky pattern} (BZ-pattern) of size $N$ is the restriction
of a BZ-function to the triangle $T_N$.

For example, a BZ-pattern of size~$4$ is an array
of nonnegative integer numbers $a_1,\dots,a_{18}$ (arranged in a triangle
as shown below) such that the numbers in any of the three hexagons satisfy 
the hexagon condition.

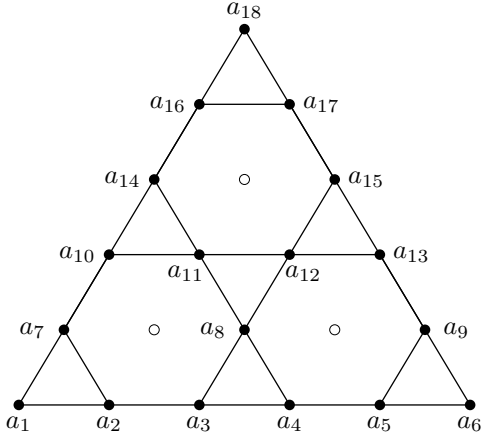
\psset{xunit=.6cm,yunit=.5cm}
\begin{center}
\pspicture[.1](-3,-1)(14,11)

\qdisk(0,0){2pt}
\qdisk(2,0){2pt}
\qdisk(4,0){2pt}
\qdisk(6,0){2pt}
\qdisk(8,0){2pt}
\qdisk(10,0){2pt}
\qdisk(1,2){2pt}
\pscircle[linewidth=.2pt](3,2){2pt}
\qdisk(5,2){2pt}
\pscircle[linewidth=.2pt](7,2){2pt}
\qdisk(9,2){2pt}
\qdisk(2,4){2pt}
\qdisk(4,4){2pt}
\qdisk(6,4){2pt}
\qdisk(8,4){2pt}
\qdisk(3,6){2pt}
\pscircle[linewidth=.2pt](5,6){2pt}
\qdisk(7,6){2pt}
\qdisk(4,8){2pt}
\qdisk(6,8){2pt}
\qdisk(5,10){2pt}

\psline[linewidth=.5pt,showpoints=false]{-}(0,0)(5,10)(10,0)(0,0)
\psline[linewidth=.5pt,showpoints=false]{-}(2,0)(1,2)(2,4)(8,4)(9,2)(8,0)
\psline[linewidth=.5pt,showpoints=false]{-}(4,0)(7,6)(6,8)(4,8)(3,6)(6,0)

\rput[u](0,-0.5){$a_1$}
\rput[u](2,-0.5){$a_2$}
\rput[u](4,-0.5){$a_3$}
\rput[u](6,-0.5){$a_4$}
\rput[u](8,-0.5){$a_5$}
\rput[u](10,-0.5){$a_6$}
\rput[u](.3,2){$a_7$}
\rput[u](4.3,2){$a_8$}
\rput[u](9.7,2){$a_9$}
\rput[u](1.3,4){$a_{10}$}
\rput[u](3.7,3.5){$a_{11}$}
\rput[u](6.3,3.5){$a_{12}$}
\rput[u](8.7,4){$a_{13}$}
\rput[u](2.3,6){$a_{14}$}
\rput[u](7.7,6){$a_{15}$}
\rput[u](3.3,8){$a_{16}$}
\rput[u](6.7,8){$a_{17}$}
\rput[u](5,10.5){$a_{18}$}
\endpspicture
\end{center}

For a BZ-pattern of size $N$, let $a_1,\dots,a_{2N-2}$ be 
the number in the lower row; $b_1,\dots,b_{2N-2}$ be 
the numbers on its left side;
and $c_1,\dots,c_{2N-2}$ be the numbers on its right side
(in all cases we count the numbers from left to right).
For the triangle on the picture above,
$b_1=a_1, b_2=a_7, b_3=a_{10},b_4=a_{14},b_5=a_{16}, b_6=a_{18}$,
and $c_1=a_{18}, c_2=a_{17}, c_3=a_{15},c_4=a_{13},c_5=a_9, c_6=a_6$.

Berenstein and Zelevinsky~\cite{BZ1} found the following interpretation of the
Littlewood-Richardson coefficients in terms of BZ-patterns.

\begin{theorem} {\rm \cite{BZ1}} \
Let $\lambda,\mu,\nu$ be three dominant weights for $GL(N)$ such that
$|\lambda|+|\mu|=|\nu|$ and
$\lambda = l_1 \omega_1 + \cdots  + l_N\omega_N$, 
$\mu = m_1 \omega_1 + \cdots +m_N\omega_N$, 
and $\nu = n_1 \omega_1 + \cdots +n_N\omega_N$.
Then the Littlewood-Richardson
coefficient $c_{\lambda\mu}^\nu$ is equal to the number of 
BZ-patterns of size $N$ with the following boundary conditions:
$$
\begin{array}{l}
l_1=b_1+b_2,\ l_2=b_3+b_4,\ \dots,\ l_{N-1}=b_{2N-3}+b_{2N-2},\\[.05in]
m_1=c_1+c_2,\ m_2=c_3+c_4,\ \dots,\ m_{N-1}=c_{2N-3}+c_{2N-2},\\[.05in]
n_1=a_1+a_2,\ n_2=a_3+a_4,\ \dots,\ n_{N-1}=a_{2N-3}+a_{2N-2}.
\end{array}
$$
\label{th:BZ}
\end{theorem}

Proposition~\ref{prop:web=BZ} says that web functions are essentially
BZ-patterns of infinite size.  
Let us fix a set of boundary rays that satisfy the conditions of 
Corollary~\ref{cor:web}.
Then an integral web functions with these boundary rays is determined by its
restriction to $T_N$, which is a BZ-pattern.  The conditions on the rays of
such web functions transform into the boundary conditions for the BZ-patterns
from Theorem~\ref{th:BZ}.
Thus Corollary~\ref{cor:web} is equivalent to Theorem~\ref{th:BZ}.

\begin{proposition}
Let $\lambda$, $\mu$, and $\nu$ be three dominant weights 
for $GL(N)$ such that $|\lambda|+|\mu|=|\nu|$.
The integral web functions that satisfy the conditions of
Corollary~\ref{cor:web} are in one-to-one correspondence with the 
BZ-patterns of size $N$ that satisfies the boundary conditions of
Theorem~\ref{th:BZ}.
This correspondence $\kappa_N$ is given by restricting a web function
$f:\RR^2\to\RR_+$ to the BZ-triangle $T_N$:
$$
\kappa_N\ :\ f \longmapsto f\,|_{T_N}\,.
$$
\end{proposition}

The following picture illustrates the statement of the proposition.
It shows a web diagram and the BZ-triangle $T_{11}$.  
In this case the corresponding BZ-pattern has $1$'s 
at the points that belong to the web diagram and it has $0$'s
everywhere else.

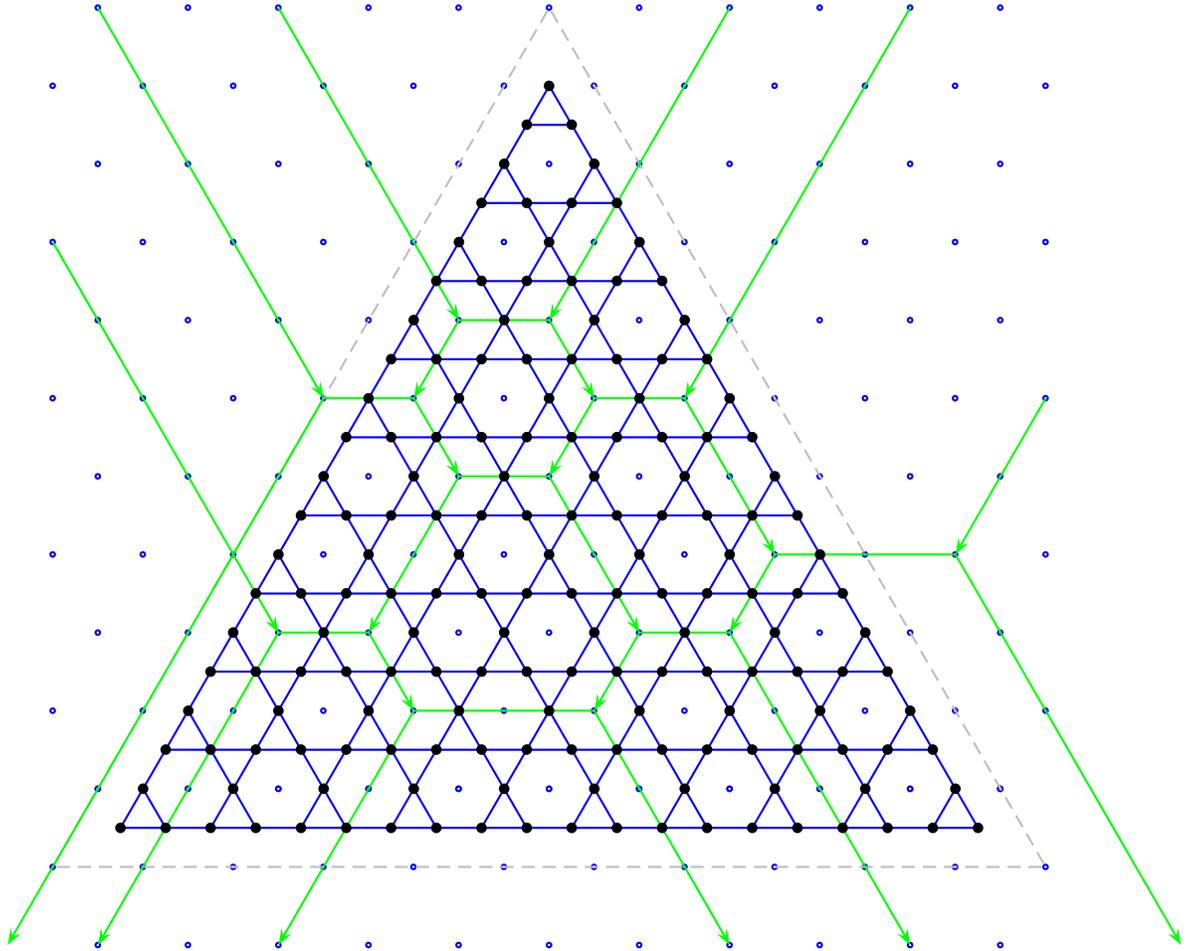
\begin{center}
\psset{unit=0.6cm,linecolor=blue,fillstyle=none,fillcolor=white}
\pspicture(-4,-4)(21,19)
\pspolygon(0,0)(19,0)(9.5,16.454)
\psline(1,1.732)(18,1.732)
\psline(2,3.464)(17,3.464)
\psline(3,5.196)(16,5.196)
\psline(4,6.928)(15,6.928)
\psline(5,8.66)(14,8.66)
\psline(6,10.392)(13,10.392)
\psline(7,12.124)(12,12.124)
\psline(8,13.856)(11,13.856)
\psline(9,15.588)(10,15.588)
\psline(1,0)(0.5,0.866)
\psline(3,0)(1.5,2.598)
\psline(5,0)(2.5,4.33)
\psline(7,0)(3.5,6.062)
\psline(9,0)(4.5,7.794)
\psline(11,0)(5.5,9.526)
\psline(13,0)(6.5,11.258)
\psline(15,0)(7.5,12.99)
\psline(17,0)(8.5,14.722)
\psline(2,0)(10.5,14.722)
\psline(4,0)(11.5,12.99)
\psline(6,0)(12.5,11.258)
\psline(8,0)(13.5,9.526)
\psline(10,0)(14.5,7.794)
\psline(12,0)(15.5,6.062)
\psline(14,0)(16.5,4.33)
\psline(16,0)(17.5,2.598)
\psline(18,0)(18.5,0.866)
\psset{border=0.1}
\multiput(0,0)(0,3.464){6}{
\multirput(0,0)(2,0){12}{\psdots[dotstyle=o](-1.7,-0.866)}}
\multiput(0,0)(0,3.464){6}{
\multirput(0,0)(2,0){11}{\psdots[dotstyle=o](-0.7,0.866)}}
\multiput(0,0)(0,-3.464){2}{
\multirput(0,0)(2,0){11}{\psdots[dotstyle=o](-0.7,0.866)}}
\psset{border=0}
\pspolygon[linestyle=dashed,linecolor=lightgray,fillstyle=none](-1.5,-0.866)
(20.5,-0.866)(9.5,18.186)
\psset{linecolor=green}
\psline(6.5,2.598)(10.5,2.598)
\psline{->}(6.5,2.598)(3.5,-2.598)
\psline{->}(10.5,2.598)(13.5,-2.598)
\psline{->}(11.5,4.33)(10.5,2.598)
\psline(11.5,4.33)(13.5,4.33)
\psline{->}(13.5,4.33)(17.5,-2.598)
\psline{->}(14.5,6.062)(13.5,4.33)
\psline(14.5,6.062)(18.5,6.062)
\psline{->}(18.5,6.062)(23.5,-2.598)
\psline{->}(20.5,9.526)(18.5,6.062)
\psline{->}(9.5,7.794)(11.5,4.33)
\psline{->}(5.5,4.33)(6.5,2.598)
\psline(5.5,4.33)(3.5,4.33)
\psline{->}(3.5,4.33)(-0.5,-2.598)
\psline{->}(7.5,7.794)(5.5,4.33)
\psline(7.5,7.794)(9.5,7.794)
\psline{->}(12.5,9.526)(14.5,6.062)
\psline{->}(10.5,9.526)(9.5,7.794)
\psline(10.5,9.526)(12.5,9.526)
\psline{->}(9.5,11.258)(10.5,9.526)
\psline{->}(6.5,9.526)(7.5,7.794)
\psline{->}(7.5,11.258)(6.5,9.526)
\psline(7.5,11.258)(9.5,11.258)
\psline(4.5,9.526)(6.5,9.526)
\psline{->}(3.5,18.186)(7.5,11.258)
\psline{->}(17.5,18.186)(12.5,9.526)
\psline{->}(13.5,18.186)(9.5,11.258)
\psline{->}(-0.5,18.186)(4.5,9.526)
\psline{->}(-1.5,12.99)(3.5,4.33)
\psline{->}(4.5,9.526)(-2.5,-2.598)

\psset{linecolor=black}

\multirput(0,0)(1,0){20}{\qdisk(0,0){2pt}}
\multirput(1,1.732)(1,0){18}{\qdisk(0,0){2pt}}
\multirput(2,3.464)(1,0){16}{\qdisk(0,0){2pt}}
\multirput(3,5.196)(1,0){14}{\qdisk(0,0){2pt}}
\multirput(4,6.928)(1,0){12}{\qdisk(0,0){2pt}}
\multirput(5,8.66)(1,0){10}{\qdisk(0,0){2pt}}
\multirput(6,10.392)(1,0){8}{\qdisk(0,0){2pt}}
\multirput(7,12.124)(1,0){6}{\qdisk(0,0){2pt}}
\multirput(8,13.856)(1,0){4}{\qdisk(0,0){2pt}}
\multirput(9,15.588)(1,0){2}{\qdisk(0,0){2pt}}

\multirput(0.5,0.866)(2,0){10}{\qdisk(0,0){2pt}}
\multirput(1.5,2.598)(2,0){9}{\qdisk(0,0){2pt}}
\multirput(2.5,4.33)(2,0){8}{\qdisk(0,0){2pt}}
\multirput(3.5,6.062)(2,0){7}{\qdisk(0,0){2pt}}
\multirput(4.5,7.794)(2,0){6}{\qdisk(0,0){2pt}}
\multirput(5.5,9.526)(2,0){5}{\qdisk(0,0){2pt}}
\multirput(6.5,11.258)(2,0){4}{\qdisk(0,0){2pt}}
\multirput(7.5,12.99)(2,0){3}{\qdisk(0,0){2pt}}
\multirput(8.5,14.722)(2,0){2}{\qdisk(0,0){2pt}}
\multirput(9.5,16.454)(2,0){1}{\qdisk(0,0){2pt}}

\endpspicture
\end{center}

\section{Remarks and Open Questions}
\label{sec:remarks}
\neweq

There are several questions that remained outside the scope of the paper.
We briefly mention them here, and they will be properly illuminated 
in subsequent publications.

First of all, an open problem of interest is to describe explicitly the transformation
maps $T_\a^{\b}$ for any two reduced words 
for~$w$ (see Section~\ref{sec:yb}).

There is an analogy between piecewise-linear transformations $T_{ijk}$ given
by~(\ref{eq:transform_general}) and the transformations for Lusztig's
parametrization of the canonical basis in $U_q^+(sl_n)$.  Lusztig's
transformations were thoroughly investigated in~\cite{BFZ}.  The combinatorial
essence of this work lies in a certain {\it chamber ansatz.}  It would be
interesting to find analogues of the results of~\cite{BFZ}.

In a forthcoming paper, Berenstein and Zelevinsky investigate string cones and
relations between Lusztig's and Kashiwara's parametrizations.  It would be
interesting to find a relationship between our combinatorial description of the
string cone in terms of rigorous paths and their construction.

Following~\cite{BFZ}, it is possible to formulate the transition maps
$T_{ijk}$ and $T_\a^\b$ in the language of the {\it tropical semiring}---%
a kind of algebraic system where one is allowed to add, multiply, and divide, but not 
subtract.  Taking the presentation of the tropical multiplication by 
the usual addition, tropical division by usual subtraction, and tropical
addition by the operation~$\min$, we can recover piecewise-linear combinatorics.
On the other hand, taking the more natural presentation of the tropical multiplication by 
the usual multiplication, tropical division by the usual division, and tropical
addition by the usual addition, we can move to the area of rational mathematics.
Hopefully, the rational expressions corresponding to the piecewise-linear transition maps
$T_\a^\b$ can be presented by some determinant-like creatures.

Knutson and Tao~\cite{Knutson} defined honeycombs as certain embeddings of
certain graphs into $\RR^2$.  They used honeycombs in the proof of Klyachko's
saturation conjecture.  Our web functions are related to honeycombs, but they
are defined in a different way by means of local conditions.  Sometimes this
definition is more convenient.  It is possible to give a proof to the
saturation conjecture in terms of web functions which is simplier than Knutson
and Tao's proof.

\medskip

\medskip

\end{document}